\documentclass[10pt]{amsart}
\usepackage{subcaption}
\usepackage{graphicx,xcolor,verbatim}
\usepackage{amssymb,amsfonts,amsmath}
\usepackage{graphics,epsfig}
\usepackage{graphicx}
\usepackage{epstopdf}
\usepackage{latexsym}
\usepackage{float}
 \usepackage{bm}

\newcommand{\be}{\begin{equation}}
\newcommand{\ee}{\end{equation}}

\newtheorem{remark}{Remark}

\newtheorem{Proposition}{Proposition}

\newtheorem{example}{Example}

\begin{document}

\title[Multi-domain chemotaxis model]{Mass-preserving approximation of a chemotaxis multi-domain transmission model for microfluidic chips}

\author{E. C. Braun} 
\thanks{Universit\'a di Roma 3, Rome, Italy, elishanchristian.braun@uniroma3.it}
\author{G. Bretti} 
 \author{R. Natalini} 
\thanks{Istituto per le Applicazioni del Calcolo ``M.Picone'', Rome, Italy, c.braun@iac.cnr.it, g.bretti@iac.cnr.it, roberto.natalini@cnr.it}
\begin{abstract}
The present work was inspired by the recent developments in laboratory experiments made on chip, where culturing of multiple cell species was possible. The model is based on coupled reaction-diffusion-transport equations with chemotaxis, and takes into account the interactions among cell populations and the possibility of drug administration for drug testing effects.\\ 
Our effort was devoted to the development of a simulation tool that is able to reproduce the chemotactic movement and the interactions between different cell species (immune and cancer cells) living in microfluidic chip environment.
The main issues faced in this work are the introduction of mass-preserving and positivity-preserving conditions involving the balancing of incoming and outgoing fluxes passing through interfaces between 2D and 1D domains of the chip and the development of mass-preserving and positivity preserving numerical conditions at the external boundaries and at the interfaces between 2D and 1D domains.
\end{abstract}



\maketitle

\section{Introduction}

The aim of the present work is to study the modelling and numerics of a chemotaxis-reaction-diffusion mathematical model describing the qualitative behavior of different cell species living in a confined environment. 
This work was inspired by laboratory experiments made on microfluidic chip \cite{V}, where some populations cohexist and interact.
In recent years, indeed, there was the development of a new approach to biological studies aimed at reconstructing organs and complex biological processes on-chip \cite{B}. The fundamental idea is that the comprehension of the sophisticated physiology of organisms, based on the complex behavior and interaction of cell populations, tissues and organs, needs interdisciplinary contributions, from biology to mathematics.

\begin{figure}[h!]
\begin{center}
\includegraphics[scale=0.5]{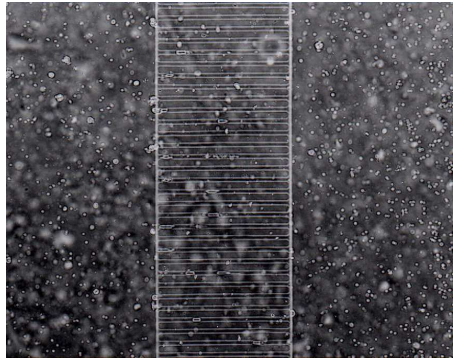}
\caption{Microfluidic chip environment: two chambers connected by multiple channels. Credits by Vacchelli et al \cite{V} edited by AAAS.}
\label{fig:chip}
\end{center}
\end{figure}

Motivated by laboratory setting of the experiment in microfluidic chips \cite{B,P,V},  we introduce a model describing the interactions  between two cells populations, namely immune and cancer cells.
 The microfluidic chip is represented as a network of channels connecting two boxes (the microfluidic chambers), see Fig. \ref{fig:chip} and a schematic picture of the experiment in Fig. \ref{fig:1}. The mathematical model, proposed in section \ref{sec:mathmod}, is a reaction-diffusion system with chemotaxis and it describes birth and death processes, migration of immune cells driven by chemical signals produced by tumor cells, interaction between different cell species.\\ 
From the mathematical point of view, we follow the framework of the classical macroscopic models of chemotaxis, where the evolution of the density of cells is described by a parabolic equation and the concentration of a chemoattractant can be given by a parabolic or elliptic equation, depending on the different regimes to be described and on authors' choices.  The choice of a continuous model to reproduce an experiment in a confined environment, with a relatively small number of cells, is motivated by the fact that we aim at developing a simulation tool which is able to describe the phenomena of immunosorveillance of cancer in tissues, where billions of cells are present. For this reason a macroscopic model is more suitable respect a particle model.\\
In the chambers we consider a 2D doubly-parabolic model which is a modification of the Keller-Segel model \cite{KS} to take into account the presence of two populations both producing chemical signal which are interacting each other. We remark that we consider only the 2D case since the experimental data do not take into account the height of the chip. Clearly, in principle, our framework could be easily extended to the third dimension.\\
 For the 1D microchannels connecting the 2D chambers we choose two different approaches: we can assign a 1D version of doubly-parabolic model used in the chambers; otherwise, we can assign a model derived from 1D-GA model \cite{GreenbergAlt87}, being  characterized by the more realistic feature that the speed of propagation of cells in the channels is finite, which seems the dominant property at this scale. 
On the other hand, other models based on hyperbolic/kinetic equations for the evolution of the density of individuals can be assigned, characterized by a finite speed of propagation \cite{GambaEtAl03,Pe, FilbetLaurencotPerthame, DolakHillen03, DiRusso}.

\subsection{Original contribution of the present paper}

From the mathematical and numerical viewpoint, here we deal with a challenging issue arising in chemotaxis modelling of cell interaction. The problem involves doubly-parabolic models in 2D domains (microfluifidic chambers) that are connected with 1D domains represented by channels, where either a doubly-parabolic or a hyperbolic-parabolic model can be assigned.  
The classical doubly-parabolic Keller-Segel (KS) model \cite{KS} of chemotaxis reads as:
 \begin{equation}\label{keller-segel}
 \left\{\begin{array}{ll}
  u_t = div \left(\nu \nabla u - \chi(u,\phi) \nabla\phi \right)\\ 
   \phi_t = D \Delta \phi + a u - b \phi ,
 \end{array}\right.
\end{equation} 
with $u$ the density of individuals in the considered medium, $\nu$ the diffusion rate of the organism according to Fick's Law and $\phi$ the density of chemoattractant. The positive constant $D$ is the diffusion coefficient of the chemoattractant; the positive coefficients $a$ and $b$, are respectively its production and degradation rates, and $\chi$ is the chemotactic sensitivity, depending on the density of the considered quantities. In the 2D domains given by the microfluidic chambers we apply a reaction-diffusion chemotaxis KS-like model inspired by (\ref{keller-segel}) and described in \ref{model}.\\ In the 1D microfluidic channels, we use the one-dimensional version of the KS-like model used in the chambers, but we also studied the behavior of individuals when a hyperbolic-parabolic model, characterized by finite speed of propagation is assigned.
Such hyperbolic-parabolic model, described in \ref{model}, is inspired by the Greeberg-Alt (GA) model, arising as a simple model for chemotaxis on a line:
\begin{equation}\label{hyper-gen}
 \left\{\begin{array}{ll}
 \partial_t u + \partial_x v =0,\\ 
  \partial_t v +\lambda^2 \partial_{x} u = -v  + \chi(u,\phi) \partial_x \phi,\\ 
\partial_t \phi = D \partial_{xx} \phi + a u - b \phi.
 \end{array}\right.
\end{equation}
 Note that here $v$ is the averaged flux.  Let us underline that the flux $v$ in model (\ref{hyper-gen}) corresponds to $v= -\lambda^2 \nabla u + \chi(u,\phi) \nabla \phi$ for the KS system. This system was analytically studied on the whole line and on bounded intervals in \cite{gumanari}, while an effective numerical approximation, the Asymptotic High Order (AHO) scheme, was introduced in \cite{NR}, see also \cite{Gosse10, Gosse11} and extended on networks with general boundary conditions in \cite{BN18} and \cite{BNR14}.
 
Here, in the numerical treatment for the computation of numerical solutions one has to take care of what happens at the inner boundaries with the switching from 2D-doubly-parabolic models and 1D-doubly-parabolic or 1D-hyperbolic-parabolic ones.

Since we aim at reproducing the numerical solutions of such models, we need to deal with a multi-domain problem given by the passage from a 2D domain represented by the chambers of the chip to 1D domains given by the channels. For this reason we need to develop ad hoc transmission conditions to ensure the mass conservation at the 2D-1D interfaces. From the numerical viewpoint, here we consider numerical boundary conditions including in the stencil a ghost cell value taken from the neighbouring domain, as we will show in the numerical Section \ref{sec:num}.
The approximation of doubly-parabolic chemotaxis models for the 1D-KS model (\ref{keller-segel}) on networks was already considered in \cite{BGKS}. However, in that case the transmission conditions were between 1D-1D interfaces and on each arc of the network the same fully-parabolic model was considered. We also underline that in such work, transmission conditions require to impose the continuity of the density of both cells $u$ and chemoattractant  $\phi$, while we  only impose the continuity of the fluxes, which seems to be more realistic when dealing with flux of individuals or molecules. \\
For the numerical approximation of the GA system (\ref{hyper-gen}), we refer to our previous papers \cite{NR} for a single line, where the numerical treatment of the hyperbolic part of the system was based on the AHO scheme with the development of mass-preserving numerical scheme at outer boundaries, while the parabolic part was approximated by finite difference and Crank-Nicolson scheme. In \cite{BN18} and \cite{BNR14} the GA system was solved on networks, thus making necessary to develop mass-preserving transmission conditions at inner nodes and suitable boundary conditions at outer nodes. However, the study of transmission conditions only involved the mass exchange between 1D-1D interfaces; moreover, on each arc of the network the same model was considered. Furthermore, the second order numerical approximation of the boundary conditions developed in such papers did not ensure the posivity preserving property in case of obscillating functions.\\

The numerical approximation of permeability Kedem-Katchalsky \cite{KK} conditions describing the conservation of the flux through a node was already considered in \cite{Q}, but we underline that in the mentioned paper the study was done for the approximation with finite elements methods of linear problems.
For reaction-diffusion problems the approximation of permeability conditions was studied in \cite{Seraf} for finite difference schemes and in \cite{Can} for discontinuous Galerkin methods. The numerical treatment of permeability conditions for chemotaxis problems was presented for the first time in \cite{DiCost2} for the 1D parabolic-parabolic interface, and a finite difference approximation was developed without taking into consideration the mass-preservation nor the positivity-preservation properties at inner nodes. 
 Therefore, to our knowledge, the present paper is the first numerical work where this new technique 
 of switching size of the domains and type of equations (parabolic vs hyperbolic approach) is introduced, in order to develop mass-preserving and positivity preserving schemes. 
 
\subsection{Main contents and plan of the paper}
In the present paper, a positivity-preserving and mass-preserving numerical discretization of Neumann boundary conditions at the corners and at the bottom and top boundaries of the 2D domain for 2D-doubly parabolic reaction-diffusion problem are presented.
Moreover, a positivity-preserving and mass-preserving numerical scheme at the inner nodes of the network connecting the 2D chambers with the 1D channels (where the 1D-doubly-parabolic or 1D-hyperbolic-parabolic problem can be assigned) is developed.
To summarize the main contents of the present work, the mathematical issues faced in this study are indentified into two aspects:
\begin{itemize}
\item the study of the behavior of two different modelling of the dynamics in the channels: the parabolic model describing the dynamics inside the chambers was coupled both with KS-like and GA-like models;
\item  the numerical approximation of equations defined in a heterogeneous domain, characterized by the switch from 2D domains, represented by microfluidic left and right chambers, to 1D domains, given by the channels connecting them.
\end{itemize}
Then, the numerical questions arising in the mentioned issues and here addressed are:
\begin{itemize}
\item the study of positivity and mass-preserving external boundary conditions for 2D-doubly-parabolic model (\ref{eqsystem});
\item the introduction of mass-preserving and positivity-preserving permeability conditions at the interfaces between 2D and 1D parabolic models, see paragraph \ref{sec:permpar};
\item the introduction of mass-preserving and positivity-preserving permeability conditions at the interfaces between 2D-fully-parabolic model and 1D-hyperbolic-parabolic model, see paragraph \ref{sec:permiper}.
\end{itemize}
        
The plan of the paper is as follows.
In section \ref{sec:bio} we describe the biological framework that inspired our study, while in section \ref{sec:mathmod} we introduce the mathematical formulation of biologically inspired models and we introduce the adopted model. Section \ref{sec:num} is devoted to the numerical techniques used to approximate the problem and in section \ref{sec:tests} some numerical tests showing the qualitative behavior of cells in the designed environment are presented. Finally, in section \ref{sec:concl} a discussion on the results and the future developments of our work is presented.

\section{Biological framework}\label{sec:bio}

The control of immune cells migration and interaction with tumor cells living inside the chambers of the microfluidic chip, represent a new and attractive approach for the clinical management of tumor deseases.   Furthermore, in the chip environment also drug testing can be exploited. Then, the quantitative assessment of immune cell migration ability to recognize and attack the tumor cells for each patient could provide a new potential parameter predictive of patient outcomes in the future.\\
Migrating cells respond to complex chemical stimuli (as mixture of growth factors, cytokines and chemokines) representing a source of chemoattractants. These chemoattractants, through the interaction with their receptors allow cells to acquire a polarized morphology and to perform the action of immunosorveillance.\\
The development of lab-on-chip technologies made it possible to realize a reproducible tailoring
of the cellular microenvironment, thus allowing the continuous monitoring of experiments
and accurate control of experimental parameters. Recently, the development of microengineering
has given the possibility to realize culturing of multiple cell types and made it possible to observe cell-cell interactions and to transpose in vivo studies to a second generation of
in vitro smart environments. The main advantages of this new technological tool are a close control over local experimental conditions and lower costs with respect to the use of animals in laboratory experiments for efficacy and toxicity testing. Some results obtained with on-chip experiments are presented in \cite{A, B, L, N,  V}.\\
Regarding the structure of microfluidic devices, they are designed to allow chemical and physical contacts between tumor cells and non-adherent immune cells (i.e. murine splenocytes or human peripheral blood mononuclear cells). The microfluidic co-culture platforms are fabricated in
polydimethylsiloxane (PDMS, Silgard 184), a biocompatible optically transparent silicone elastomer.\\
In particular, here we refer to the experiment of two main culture chambers (a tumor and an immune cell compartment) connected via narrow capillary migration micro-channels having, respectively, width, length and height of $12 \mu m$, $500 \mu m$ and $10 \mu m$.  The cross-sectional dimensions of culture compartments are $1 mm$ (width) $\times$ $100 \mu m$ (height).
 For the development of the mathematical modelling explained in the next section we remark that we neglect the third dimension, thus we consider corridors and chambers as 2D objects.
  Two populations, immune and cancer cells, are introduced in two separate chambers; in particular, the immune cells are in the right chamber and the cancer cells in the left one (target chamber). The microchannels connect the two areas and allow the chemical diffusion and the migration. 
  The culture medium is neutral, thus meaning that no exogenous substance are introduced. Mainly, the dynamics observed is the migration of immune cells from the right to the left in order to attack the tumor cells. The laboratory experiments are made in the context of immune competence vs. immunodeficiency, i.e. in a healthy or
defective immune system. Indeed, the complex interactions, including cell-cell contacts, between cancer cells and immune system, which acts by limiting or suppressing tumor progression, is crucial in the tumor growth and invasion process.\\

\begin{figure}
\begin{center}
$$\input{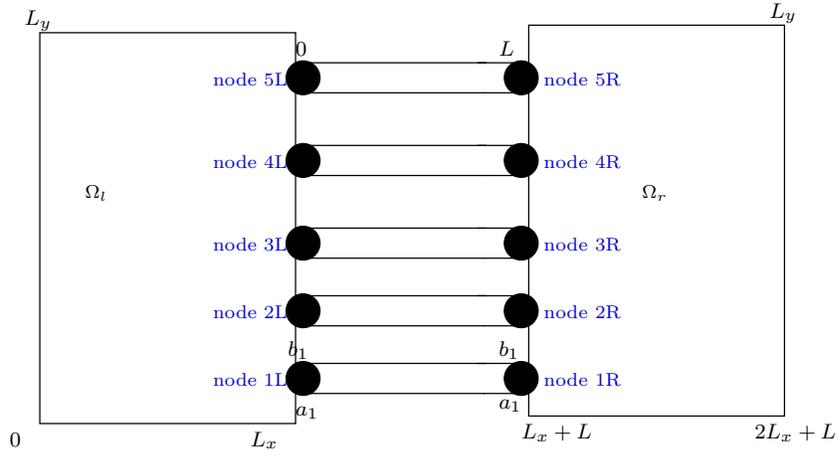}$$
\end{center}
\caption{Simplified schematization of the chip geometry depicted in Fig. \ref{fig:chip}.}
\label{fig:1}
\end{figure}

\section{Mathematical framework} \label{sec:mathmod}
Nowadays, mathematical analysis of biological phenomena has become an important tool to explore
complex processes, and to detect mechanisms that might not be evident to the experimenters.
Although a mathematical model cannot replace a real experiment, it may 
represent a support tool to explain acquired biological data and it may consent to gain a deeper understanding of the interactions between cancer cells and immune system.
More generally, mathematical models can describe a broad variety of biological phenomena, including cell dynamics and cancer \cite{12}, \cite{16}, \cite{13}, \cite{14}, \cite{15}.\\
The movement of bacteria under the effect of a chemical substance has been widely studied in the last decades, and numerous mathematical models have been proposed. As shown in \cite{murray}, chemotaxis is decisive in biological processes. For instance, the formation of cells aggregations (amoebae, bacteria, etc) occurs during the response of the different species to the change of the chemical gradients in the environment. 
 Moreover it is possible to describe this biological phenomenon at different scales. For example, by considering the population density as a whole, it is possible to obtain macroscopic models of partial differential equations.\\ 

In this paper, in order to describe the dynamics of cells in the 2D chambers we use a KS-like model, while in the microchannels we compare the behavior between two different modelization: 1D KS-like model and 1D GA-like model. The modelling here applied is described in the next subsection \ref{model}.
\subsection{The model}\label{model}

Here we introduce a mathematical model that aims at describing the behavior of two populations of cells cohexisting together: tumoral cells $T$ and immune cells(macrophage) $M$. We underline that the setting here considered can be make more complex with the introduction of a greater number of cell species and with the presence of an exogenous substance in the environment.\\ 
The model consists of a reaction-diffusion system with chemotaxis, that it is able to describe birth and death processes, interaction with chemoattractant, interaction and competition between different cell species.  The microfluidic chip is schematized as a network of channels connecting two boxes (the microfluidic chambers), then, following the ideas in \cite{BN18}, ad hoc transmission conditions were introduced to ensure the mass conservation. The parameters of the model, such as the velocity of different cell populations, the turning rates and the decay rates will be calibrated with observed data.

Cancer cells $T$ produce chemical signal, called $\varphi$, activating the immune response of $M$  and influencing their behavior.
Moreover, we take into account the presence of cytokines $\omega$ (produced by $M$), acting as a chemical killer of cancer cells. 
Then, the model describing the dynamics of the two cell species and the diffusion of the chemoattractant in the 2D chambers reads as:
 
\begin{equation}\label{eqsystem}
 \left\{\begin{array}{ll}
\frac{\partial}{\partial t}  T = D_T \Delta T - \lambda_{T}(\omega)T - k_T(t) T,\\
\\
\frac{\partial}{\partial t} M = D_M  \Delta M - div(\chi(M,\varphi)\nabla \varphi) - k_M(t) M,\\
\\
\frac{\partial}{\partial t} \varphi = D_{\varphi}  \Delta \varphi + \alpha_\phi T -\beta_{\varphi} \varphi,\\
\\
\frac{\partial}{\partial t} \omega =  D_{\omega}  \Delta \omega + \alpha_\omega M -\beta_\omega \omega,\\
\\
 \end{array}\right.
\end{equation}
and we need to assign suitable initial conditions and boundary conditions for the cells and the chemoattractant concentrations that will be specified in the next paragraphs.\\ 
In particular, the system above describes the following situation: tumor cells $T$ produce a chemical substance $\varphi$ attracting immune cells $M$ and enabling them to recognize and interact with tumor cells. Immune cells also produce a chemical substance $\omega$ which makes the immune cells able to migrate towards the tumor cells. Therefore, in the first equation of the system (\ref{eqsystem}), besides the diffusion term, we can find $-\lambda_{T}(\omega)T$ representing the tumor suppression operated by immune cells. In the second equation, in addition to the diffusion term we have the chemotactic term $f= \chi(M,\varphi)\nabla \varphi$ due to the presence of the chemical substance $\varphi$ produced by the tumor. We remark that both in the first and in the second equation we include a term $- k_T(t) T$ and $-k_M(t) u$ taking into account the possibility of drug administration, with the functions $k_T$ and $k_M$ having an exponential decay in time:
\begin{equation}
\begin{split}
 k_T (t) &=  K_T e^{-\alpha_T t},\\
k_M (t) &=  K_M e^{-\alpha_M t}.
\end{split}
\end{equation}

We also need to introduce the functions:
\begin{equation}
\begin{split}
 \chi (M,\varphi) &=  \frac{k_1 M}{(k_2 + \varphi)^\gamma},\\
\lambda_{T}(\omega) &= \frac{k_\omega \omega}{1 + \omega},\\
\end{split}
\end{equation}
representing, respectively, the chemotactic sensitivity of immune cells and the decay rate of cancer cells under the action of immune cells. Note that $k_\omega$ is the killing efficiency of immune cells, $k_1$ represents the cellular drift velocity, while $k_2$ is the receptor dissociation constant, which says how many molecules are necessary to bind the receptors. We mainly refer to \cite{murray} for the values of the parameters $k_1$, $k_2$, $\gamma$, and all the parameters are reported in Table \ref{table:param1}. \\
Now, in order to describe the dynamics of cells in the microchannels connecting the two boxes, we introduce the following 1D models for the dynamics. To this aim, we consider two possible approaches to observe a different dynamics in the channels:
\begin{itemize}
\item if we assign the 1D doubly-parabolic model on each channel, we have one-dimensional version of system (\ref{eqsystem}), where the superscript $c$ indicates that we consider all the quantitites in the channels:
\begin{equation}\label{eqsystem1D}
 \left\{\begin{array}{ll}
\frac{\partial}{\partial t}  T_c = D_T \Delta T_c - \lambda_{T}(\omega)T_c - k_T(t) T_c,\\
\\
\frac{\partial}{\partial t} M_c = D_M  \Delta M_c - \partial_x f_c - k_M(t) M_c,\\
\\
\frac{\partial}{\partial t} \varphi_c = D_{\varphi}  \Delta \varphi_c + \alpha_\varphi T_c -\beta_{\varphi} \varphi_c,\\
\\
\frac{\partial}{\partial t} \omega_c =  D_{\omega}  \Delta \omega_c + \alpha_\omega M_c -\beta_\omega \omega_c.\\
\\
 \end{array}\right.
\end{equation}
\item if we consider the 1D hyperbolic-parabolic model on each channel, we have the following system: 

\begin{equation}\label{GA1D}
 \left\{\begin{array}{ll}
 \partial_t T_c + \partial_x v_c^T = -\lambda_T(\omega) T_c - k_T(t) T_c,\\
\partial_t v_c^T + D_T \partial_x T_c =  - v_c^T,\\
\partial_t \omega_c = D_{\omega_c}  \partial_{xx} \omega_c + \alpha_\omega T_c -\beta_c \omega_c,\\
 \partial_x M_c + \partial_t v_c^M = -k_M(t) M_c,\\
\partial_t v_c^M + D_{M_c} \partial_x M_c = f_c - v_c^M,\\
\partial_t\varphi_c = D_{\varphi_c}  \partial_{xx}\varphi_c + \alpha_\varphi T_c -\beta_\varphi \varphi_c,\\
 \end{array}\right.
\end{equation} 
where $f_c = \chi(M_c,\varphi_c)\partial_x \varphi_c $ and with $v_c^T$ and $v_c^M$, respectively, the average flux of tumor cells $T_c$ and  immune cells $M_c$ in the channels. 
\end{itemize}

We remark that, for the hyperbolic-parabolic system (\ref{GA1D}) we also need to assign initial and boundary conditions for the flux $v$.

For the sake of simplicity, we write the 2D model (\ref{eqsystem}) as a general 2D-doubly-parabolic system with source term as:
   \begin{eqnarray}\label{parabolic2d}
   \left\lbrace
   \begin{array}{lcl}
   \partial_{t}u &=&D_{u}\Delta u -\text{div}f + g(x,y,t,u)\\
   \partial_{t}\phi &=&D_{\phi}\Delta \phi + a u - b \phi,\\
   \end{array}\right.
   \end{eqnarray}
   with $u$ the density of individuals and $\phi$ the density of chemoattractant. From now on, the two components of the drift term $f$ will be indicated as:
   \begin{eqnarray}
   f\left(x,y,t\right):=\left(
   \begin{array}{c}
   f^x\left(x,y,t\right)\\
   f^y\left(x,y,t\right)
   \end{array}\right).
   \end{eqnarray} 
  In the mono-dimensional channel we rewrite the 1D-doubly-parabolic system (\ref{eqsystem1D}) in a more general form:
   \begin{eqnarray}\label{parabolic1d}
   \left\lbrace\begin{array}{lcl}
   \partial_{t}u_{c} &=& D_{u_c}\partial_{xx}u - \partial_x f_{c} + g(x,t,u),\\
   \partial_{t}\phi_{c} &=&D_{\phi_{c}}\partial_{xx}\phi_{c} +a_{c} u_{c}-b_{c}\phi_c\\
   \end{array}\right.
   \end{eqnarray}
  and the 1D hyperbolic-parabolic system (\ref{GA1D}) with source term rewrites as:
   \begin{eqnarray}\label{hyperbolic1d}
   \left\lbrace\begin{array}{lcl}
  \partial_t u_c + \partial_x v_c &=& g(x,t,u_c),\\ 
  \partial_t v_c +\lambda_c^2 \partial_{x} u_c &=& -v_c  + f_c,\\ 
\partial_t \phi_c &=& D_{\phi_c} \partial_{xx} \phi_c + a_c u_c - b_c \phi_c,
 \end{array}\right.
\end{eqnarray}
  with $f_c = \chi(u_c,\varphi_c) \partial_x \phi_c$. In Table \ref{table:param1} are reported the parameters of the problem.\\
The systems above have to be complemented with initial conditions for the unknowns $u, v, \phi$, assumed to be smooth; initial data will be specified in paragraph \ref{sec:algo}. On the boundary we consider for all the quantities homogeneous Neumann conditions, so that we are assuming no-flux boundary conditions.\\

\textbf{Monotonicity conditions.}
We also mention at this point that this model has an analytical monotonicity criteria. For linear convection term $f_c=a u$, and linear source term $g=b u$ the criteria  
$$
\left\vert\frac{a}{\lambda}\right\vert -b \leq 1,
$$
must be satisfied in order for the quantity $u$ to be non-negative.
Otherwise we would have negative $u$ which would lead to unphysical solutions.\\
In regards to our model, that would mean we have for the immune cell density $M$ the monotonicity condition 
\begin{equation}
\frac{k_{1}}{\left(k_{2}+\varphi\right)^{\gamma}} \vert \partial_{x}\varphi \vert\leq\sqrt{D_{M}}
\end{equation}
 and for tumor cell density $T$:
 \begin{equation}
 \frac{k_{\omega} \omega}{1+\omega} T \leq 1 
\end{equation}
to be verified in the computational domain in order to ensure non-negative solutions.

\begin{remark}
We remark that the no-flux conditions boundary conditions used in our simulations are needed to have the mass-conservation of all the quantities. However, they are not realistic, since in the laboratory experiment there is an inflow of cells from the outer boundaries. In our future developments we will extend the no-flux boundary conditions to more general ones.
\end{remark}
For our implementation we want to model this PDE in two domains.
A simplified schematization of the bounded surface where experiment is performed is reported in Fig. \ref{fig:1}. 
We have two microfluidic chambers of the same size, one on the left and the other on the right, defined, respectively, as $\Omega_l=[0,L_x]\times[0, L_y]$ and $\Omega_{r}:=[L_{x}+L,2L_{x}+L]\times[0,L_{y}]$ they are connected by  microchannels, each of them schematized for simplicity as a line $I=[0, L]$. Thus, the link between the box on the left and the corridor is schematized as a junction (node $1L$) and analogously the link between the corridor and the box on the right, as node $2L$. The two junctions are not really a single point, thus they are parametrized as an interval for node $1L$ and node $2L$, namely $[a_1,b_1]$ of length $\sigma:=b_1-a_1$. 
We remark that for the sake of simplicity, the numerical treatment is developed for a simplest geometry composed by 2D chambers connected through a single 1D channel. The extension to multiple 1D channels is done in paragraph \ref{multi}.

\subsection{Outer and inner boundary conditions for the models with source term $g=0$.}
From now on, in order to use the mass conservation argument at the outer boundaries of the 2D domain and at the inner interface between 2D and 1D domains, we study a simplified version of the models (\ref{parabolic2d})-(\ref{parabolic1d}) and (\ref{parabolic2d})-(\ref{hyperbolic1d}) putting the source term $g$ equal to zero:
\begin{align}\label{parabolic2d1d_f0}
  a)   \left\lbrace 
 \begin{array}{lcl} 
  \partial_{t}u = D_{u}\Delta u -\text{div}f \\ 
   \partial_{t}\phi = D_{\phi}\Delta \phi + a u - b \phi,\\  
   \end{array}\right. 
   & b) \ \left\lbrace\begin{array}{lcl}
   \partial_{t}u_{c} = D_{u_c}\partial_{xx}u - \partial_x f_{c},\\
 \partial_{t}\phi_{c} = D_{\phi_{c}}\partial_{xx}\phi_{c} +a_{c} u_{c}-b_{c}\phi_c,\\
   \end{array}\right.
 \end{align}
  \begin{align}\label{parabolic2dhyper1d_f0}
    a) \left\lbrace 
    \begin{array}{lcl} 
 \partial_{t}u = D_{u}\Delta u -\text{div}f \\ 
   \partial_{t}\phi = D_{\phi}\Delta \phi + a u - b \phi,\\  
   \end{array}\right. 
   & b) \ \left\lbrace\begin{array}{lcl}
  \partial_t u_c + \partial_x v_c = 0,\\ 
  \partial_t v_c +\lambda_c^2 \partial_{x} u_c = -v_c  + f_c,\\ 
\partial_t \phi_c = D_{\phi_c} \partial_{xx} \phi_c + a_c u_c - b_c \phi_c.
 \end{array}\right.
 \end{align}
In particular, we have to prescribe the flux conservation at the inner boundaries for  the 2D-1D parabolic case (\ref{parabolic2d1d_f0}) and for the 2D parabolic-1D hyperbolic case (\ref{parabolic2dhyper1d_f0}), since we cannot loose nor gain any cells during the passage through a node.
While keeping most of the  boundary conditions, we must change them at the interface node.\\
In the 2D left box $\Omega_l$, the position of node $1L$ is at $x=L_{y}, y \in [a_1,b_1]$ and for the 1D domain represented by the channel node $1L$ is placed at $x=0$, see Fig. \ref{fig:1}.\\

   We shall note at this point that although we work in the following with the general systems (\ref{parabolic2d}),(\ref{parabolic1d}) and (\ref{hyperbolic1d}), the same results can be applied to the complete models (\ref{eqsystem}),(\ref{eqsystem1D}) and (\ref{GA1D}).

\subsubsection{Boundary conditions for the 2D doubly-parabolic model (\ref{parabolic2d1d_f0})-a)}
Considering that our model describes the migration of cells by both diffusion and chemoattractant effects, physically speaking the mass of cells and the chemoattractant must be preserved in absence of creation and destruction of cells.
   For $\phi$ in (\ref{parabolic2d1d_f0})-a), by using the divergence theorem, we can write:
   \begin{eqnarray}\label{BCchemo_analytical}
   \begin{array}{llcl}
   &\displaystyle\frac{d}{dt}\displaystyle\int_{\Omega}\phi\left(x,y,t\right)d\Omega_l=\displaystyle\int_{\Omega}D_{\phi}\triangle\phi\left(x,y,t\right)d\Omega\\
  =&\displaystyle \oint_{\delta\Omega}D_{\phi} \triangledown\phi\left(x,y,t\right)\bm{n}dS = 0.
   \end{array}
   \end{eqnarray}
    We assume no-flux condition for the chemoattractant in order to preserve its mass in absence of source terms.\\
    Since we defined the source term $f$ as a product function of $\triangledown \phi$, we get the equivalent condition $f\left(t,x,y\right)\bm{n}\vert_{\delta\Omega}=0$. Note that for $u$ the same condition holds:
\begin{eqnarray}\label{BCcells_analytical} 
  \triangledown u\left(t,x,y\right)\bm{n}\left.\right\vert_{\delta \Omega}=0, (x,y)\in \delta\Omega,
   \end{eqnarray}
and these boundary conditions guarantee mass-conservation.\\

   The same approach gives no-flux condition for both $u$ and $\phi$ in the right chamber $ \Omega_r$ (and for the complete model for $T$ and $\omega$ as well).\\
   
   
\subsubsection{Interface between 2D-1D models in (\ref{parabolic2d1d_f0})}\label{sec:permpar}
Here we prescribe the conservation of the mass between the left box and the corridor (node $1L$ in Fig. \ref{fig:1}). The conservation condition reads as:
$$
\displaystyle\frac{d}{dt}\displaystyle\int_{\Omega}u\left(x,y,t\right)d\Omega+\displaystyle\frac{d}{dt}\displaystyle\int_{0}^L u_{c}\left(x,t\right)dx =0,
$$
and it rewrites as:
\begin{eqnarray*}
\begin{array}{llcl}
 &0= \displaystyle\int_{\Omega}( D_{u}\triangle u\left(x,y,t\right)-\text{div}f \left(x,y,t\right)) d\Omega_l\ +\displaystyle\int_{0}^L (D_{u_c}\partial_{xx} u\left(x,t\right)-\partial_{x}f_{c}\left(x,t\right)) dx\\
 =& \displaystyle\oint_{\delta\Omega}\left(D_{u}\triangledown u\left(x,y,t\right)-f\left(x,y,t\right)\right)\bm{n}dS +
 \displaystyle\int_{0}^{L}\left(D_{u_c}\partial_{xx}u_{c}\left(x,t\right)-\partial_{x}f_{c}\left(x,t\right)\right)dx,
\end{array}
\end{eqnarray*}
by using the divergence theorem in the first integral.
With our analytical boundary conditions (\ref{BCcells_analytical}), the integral vanishes except at the boundary where the node is positioned.\\
We remark that attention has to be paid with $\bm{n}$ being the outer normal of the domain. We have:
\begin{eqnarray*}
\begin{array}{llcl}
& \displaystyle\int_{a_1}^{b_1}\left(D_{u}\partial_{x} u \left(L_{x},y,t\right)- f^x(L_{x},y,t)\right)dy=
-\displaystyle\int_{0}^{L}\left(D_{u_c}\partial_{xx} u_{c}\left(x,t\right)-\partial_{x}f_{c}\left(x,t\right)\right)dx,
\end{array}
\end{eqnarray*}
and, thanks to the boundary conditions (\ref{BCchemo_analytical}) and (\ref{BCcells_analytical}) some terms cancel in equation above, thus we get the condition:
 \begin{equation}
 \label{KKintegral_2D1Dparabolic}
\displaystyle\int_{a_1}^{b_1} \left(D_{u} \partial_{x} u \left(L_{x},y,t\right)-f^x(L_{x},y,t\right)) dy = D_{u_c}\partial_{x}u_{c}\left(0,t\right)-f_{c}\left(0,t\right).
 \end{equation}

Now we impose Kedem-Katchalsky (KK) \cite{KK} conditions describing the conservation of the flux through a node (see also \cite{Q} for numerical treatment of these conditions). In particular, at the interface  between left chamber and channels we have (on the left of node $1L$ in Fig. \ref{fig:1}):
\begin{equation}\label{BCKK1_2D1Dparabolic}
D_{u} \partial_{x} u \left(L_{x},y,t\right)- f^x \left(L_{x},y,t\right)=K\left(u_{c}\left(0,t\right)-u(L_{x},y,t)\right)\qquad \text{for $y\in\left[a_1,b_1\right]$}
\end{equation}
and on the right of node $1L$ we have:
\begin{eqnarray}\label{BCKK2_2D1Dparabolic}
D_{u_c}\partial_{x}u_{c}\left(0,t\right)-f_{c}\left(0,t\right)&=&\displaystyle K(u_c (0,t)) \sigma -\int_{a_1}^{b_1} u (L_{x},y,t) dy.
\end{eqnarray}
Thanks to conditions (\ref{BCKK1_2D1Dparabolic}) and (\ref{BCKK2_2D1Dparabolic}) we are guaranteed to have the flux conservation (\ref{KKintegral_2D1Dparabolic}); and we will use such conditions to obtain numerical boundary conditions for the boundary values at the nodes on both sides, as shown in Section \ref{sec:num} in paragraph \ref{sez:discrtrans}.

\subsubsection{Interface between 2D-1D models in (\ref{parabolic2dhyper1d_f0})}\label{sec:permiper}

In this section we describe the combination of 2D parabolic-1D hyperbolic model in order to describe the dynamics with a hyperbolic model (\ref{hyperbolic1d}) in the one-dimensional domain represented by microchannels. Further care has to be made in order to keep some important properties which ensure consistency and non-negativity of numerical solutions when connecting both models.\\

Now the transmission condition for the switch from $\Omega_l$ to $I=[0,L]$ are derived in this case. For the mass conservation we impose the condition:
\begin{eqnarray*}
\begin{array}{llcl}
&0 =\displaystyle\frac{d}{dt}\displaystyle\int_{\Omega}u \left(x,y,t\right)d\Omega_l+ \displaystyle\frac{d}{dt}\displaystyle\int_{0}^L u_{c}\left(x,t\right)dx\\
&= \displaystyle \int_{\Omega}\left(D_{u}\triangle u \left(x,y,t\right)-\text{div} f \left(x,y,t\right)\right)d\Omega+ \displaystyle\int_{0}^{L}-\partial_{x} v\left(x,t\right)dx\\
\Longrightarrow & \displaystyle\oint_{\delta\Omega_l}\left(D_{u} \triangledown u \left(x,y,t\right)-f(x,y,t)\bm{n}\right)dS + v\left(0,t\right)=0.\\
\end{array}
\end{eqnarray*}
 Note that in the above formula we have $v\left(L,t\right)=0$ because we are looking at left interface (node $1L$). Then, we finally get:
\begin{equation}\label{cond}
\displaystyle\int_{a_{1}}^{b_{1}} \left(D_{u} \partial_{x} u \left(L_{x},y,t\right)- f^x \left(L_{x},y,t\right)\right)dy  = -v\left(0,t\right).
\end{equation}
Now we impose the KK-condition at the interface:
\begin{equation*}
D_{u}\partial_{x} u \left(L_{x},y,t\right)- f^x \left(L_{x},y,t\right) = K\left(u_{c}\left(0,t\right)- u\left(L_{x},y,t\right)\right)\qquad \textrm{ for } y\in\left[a_{1},b_{1}\right]
\end{equation*}
and then (\ref{cond}) reads as:
\begin{eqnarray}\label{BCKK2_2D1D_hyperbolic}
v\left(0,t\right)&=&-K\sigma u_{c}\left(0,t\right)+K\displaystyle\int_{a_{1}}^{b_{1}}u_{l}\left(L_{x},y,t\right)dy.
\end{eqnarray}

\section{Numerical approximation} \label{sec:num}
Here we describe the numerical approximation of the adopted models, 2D-doubly-parabolic, 1D-doubly-parabolic and 1D-hyperbolic-parabolic.
We define equispaced $x_{i}:=i\triangle x$, $t_{n}:=n \triangle t$ and $y_{j}:=j\triangle y$ with $\triangle x$, $\triangle y$, $\triangle t>0$ and $i=0,\dots,N_{x}+1$, $j=0,\dots,N_{y}+1$; for the channel $[0,L]$ we discretize it as $x_i=i \triangle x$, with $i=0,\ldots,N$.
For a more structured presentation, we introduce the operators
\begin{eqnarray*}
\delta_{x}^{2}u_{i,j}^{n}&:=&u_{i+1,j}^{n}-2u_{i,j}^{n}+u_{i-1,j}^{n}, \ \delta_{y}^{2}u_{i,j}^{n}:=u_{i,j+1}^{n}-2u_{i,j}^{n}+u_{i,j-1}^{n},\\
\delta_{x}^{0}u_{i,j}^{n}&:=&u_{i+1,j}^{n}-u_{i-1,j}^{n}, \ \ \delta_{y}^{0}u_{i,j}^{n}:=u_{i,j+1}^{n}-u_{i,j-1}^{n},\\
\delta_{x}^{1}u_{i,j}^{n}&:=&u_{i+1,j}^{n}-u_{i,j}^{n}, \ \ \delta_{y}^{1}u_{i,j}^{n}:=u_{i,j+1}^{n}-u_{i,j}^{n}.
\end{eqnarray*}
 We use a first order explicit finite difference method in time and a second order central method for the approximation of the diffusion term in space. For the chemoattractant term we use a finite difference scheme in space as will be specified in the sequel.\\
 We remark that using a purely explicit methods implies restrictions on the mesh grid spacing and time step size to ensure stability due to the P\'eclet number criterion. However, to prevent strong restrictions on the mesh grid in the case of dominant advection regime, we introduce artificial viscosity, which leads to less restrictive condition when diffusion is small, but it decreases the order of the scheme. In the sequel we always assume to have the mesh grid small enough, thus neglecting the artificial viscosity term, which will be addressed in the numerical approximation of the model in section \ref{sec:num}.

\begin{remark}
Note that special attention has to be paid also to the source term $g\left(x,y,t,u\right)$. Although in a simple explicit method one can evaluate the function at each time step $n$ at mesh grid point $(i,j)$, the function itself can induce stiffness, enforcing small time steps. To overcome this issue, implicit methods can be used such as the Crank-Nicolson-method. However, here we work with pure explicit methods to make an easier presentation of the schemes and we will address the implicit method for these models in the sequel, see paragraph \ref{sec:algo}.
\end{remark}
 
\begin{figure}[h!]
\includegraphics[scale=0.3]{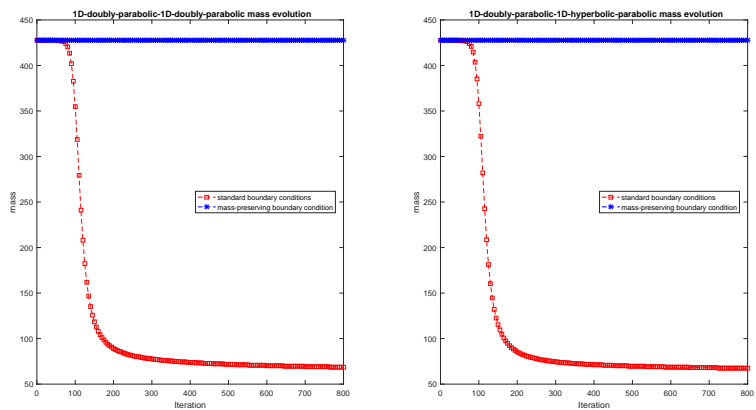}
\caption{On the left: evolution of total mass for 1D-doubly-parabolic model with standard(finite difference) and mass-preserving boundary conditions. On the right: evolution of total mass for 1D-hyperbolic-parabolic model with standard (finite difference) and mass-preserving boundary conditions.}
\label{fig:masspreserving}
\end{figure}

 Another issue is the choice of the right boundary conditions which should reflect the qualitative attributes of the analytical model.
 In absence of source terms, the mass of cells and chemical substances are preserved. In order to make a numerical verification of this property, we considered the numerical approximation at the interface between 1D-1D models. In more detail, choosing standard boundary conditions by simply discretizing Neumann boundary conditions with a finite difference scheme, the mass will not be preserved over time, see Fig. \ref{fig:masspreserving}. In particular, in Fig. \ref{fig:masspreserving} a comparison between mass-preserving and usual finite difference boundary condition is performed, for the 1D-doubly parabolic case on both sides of the interface (on the left) and for the 1D-doubly-parabolic-1D-hyperbolic-parabolic interface. From this 1D numerical example it is evident the necessity to develop modified boundary conditions which are consistent and preserve the mass correctly.\\

Mass-preserving and positivity-preserving numerical approximation will be developed in the present section.
 In the following we will neglect the label $c$ to make the reading easier and make distinction only when necessary.

\subsection{The parabolic-parabolic case}
Here we propose a numerical scheme for the doubly-parabolic systems (\ref{parabolic2d}) and (\ref{parabolic1d}).

For the discretization of equations in 2D system (\ref{parabolic2d}) in the interior points of the domain, i.e.  for $i=1,\ldots,N_x, j=1,\ldots,N_y$, we define an explicit in time finite difference discretization both for $u$ and $\phi$:
   \begin{eqnarray}\label{uparabolic_2d_discr}
   \begin{array}{lcl}
   \frac{u_{i,j}^{n+1}-u_{i,j}^{n}}{\triangle t}&=&
   \frac{D_{u}}{\triangle x^{2}}\left(u_{i-1,j}^{n}-2u_{i,j}^{n}+u_{i+1,j}^{n}\right)+\frac{D_{u}}{\triangle y^{2}}\left(u_{i,j-1}^{n}-2u_{i,j}^{n}+u_{i,j+1}^{n}\right)\nonumber\\
   &&-\triangle_{i}^{n}\left(f^x_j\right)-\triangle_{j}^{n}\left(f^y_i\right)
   \end{array}
    \end{eqnarray}
    with $\triangle_{i}^{n}\left(f^x_j\right)$ consistent approximation for $\text{div}\left(f\right)$, i.e. a second order central in space finite difference:
     \begin{eqnarray} \label{star2D}
     \Delta_{i}^{n}(f^x_j):=\frac{f^{x,n}_{i+1,j}-f^{x,n}_{i-1,j}}{2\triangle x}, \ i=1,\dots,N_x,\j=1,\ldots,N_y,
     \end{eqnarray}
    and analogously for $\triangle_{j}^{n}\left(f^y_i\right)$. Note that the function $f^x=\chi(u)\phi_{x}$ can be discretized with $f^{x,n}_{i}=\chi(u^n_{i,j})(\phi_x^n)_{i,j}$ with $(\phi_x^n)_{i,j}$ an appropriate second order approximation of $\phi_{x}$:
     \begin{equation}\label{phix}
     (\phi^n_x)_{i,j} = \frac{\phi^n_{i+1,j}-\phi^n_{i-1,j}}{2 \Delta x}, \  i=1,\dots,N_x, \ j=1,\ldots,N_y.
     \end{equation}
    For the chemoattractant we have the approximation scheme:
     \begin{equation}\label{phi2despl}
     \begin{array}{lcl}
   \frac{\phi_{i,j}^{n+1}-\phi_{i,j}^{n}}{\triangle t}&=&\frac{D_{\phi}}{\triangle x^{2}}\left(\phi_{i-1,j}^{n}-2\phi_{i,j}^{n}+\phi_{i+1,j}^{n}\right)+\frac{D_{\phi}}{\triangle y^{2}}\left(\phi_{i,j-1}^{n}-2\phi_{i,j}^{n}+\phi_{i,j+1}^{n}\right)\\
   &&+a u_{i,j}^{n}-b \phi_{i,j}^{n}.
   \end{array}
   \end{equation}
For 1D system (\ref{parabolic1d}) in the interior points of the channel we apply the same explicit in time finite difference scheme both for $u$ and $\phi$: 
   \begin{equation}\label{uparabolic1d_discr}
   \frac{u_{i}^{n+1}-u_{i}^{n}}{\triangle t} = \frac{D_{u_c}}{\triangle x^{2}}\left(u_{i-1}^{n}-2u_{i}^{n}+u_{i+1}^{n}\right)-\triangle_{i}^{n}\left(f\right), \quad i=1,\ldots,N \\
   \end{equation}
   and
    \begin{equation}\label{phiparabolic1d_discr}
   \frac{\phi_{i}^{n+1}-\phi_{i}^{n}}{\triangle t} = \frac{D_{\phi_{c}}}{\triangle x^{2}}\left(\phi_{i-1}^{n}-2\phi_{i}^{n}+\phi_{i+1}^{n}\right)+a_{c}u_{i}^{n}-b_{c}\phi_{i}^{n}, \qquad i=1,\ldots,N,
   \end{equation}
    with $\triangle_{i}^{n}\left(f\right)$, as above, a second order central in space finite difference:
     \begin{equation} \label{star}
     \Delta_{i}^{n}(f):=\frac{f^n_{i+1}-f^n_{i-1}}{2\triangle x}, \ i=1,\dots,N_x,
     \end{equation}
   where the 1D version of (\ref{phix}) for the approximation of $\phi_x$ in $f$ is used.\\

\textbf{CFL condition.}
By using the Von-Neumann stability analysis we obtain the following stability criterias (CFL-condition).
\begin{eqnarray}
D_{u_c}\frac{\triangle t}{\triangle x^{2}}\leq\frac{1}{2}\qquad\text{ for 1D},\\
D_{u}\frac{\triangle t}{\triangle x^{2}}+D_{u}\frac{\triangle t}{\triangle y^{2}}\leq\frac{1}{2}\qquad\text{ for 2D}.
\end{eqnarray}
These restrictions for the step size and the mesh grid size can be avoided by using implicit methods, but this would increase computational cost because of the necessity of solving a non-linear equation system at each iteration.

   \begin{remark}
   Since we are dealing with an explicit method, we can calculate the values for the next time step $n+1$ by using solely the values of the previous time step $n$;  we remark that an implicit approximation can be applied, with the use of Crank-Nicolson (CN) method, in order to increase the accuracy to second order and avoid restriction of $\triangle t$ and $\triangle x$ due to the CFL-Condition which do not arise with the CN-method. Indeed, we remark that the implemented algorithm for simulations described in paragraph \ref{sec:algo} is based on the CN-method in time.
   \end{remark}
  
For the two-dimensional system here we report the numerical scheme.  
  
The numerical method in the interior points of the 2D domain for the cell density reads as:
\begin{eqnarray}\label{parabolic2d_discr_complete}
\ u_{i,j}^{n+1} &=& u_{i,j}^{n}+D_{u}\frac{\triangle t}{\triangle x^{2}}\left(u_{i-1,j}^{n}-2u_{i,j}^{n}+u_{i+1,j}^{n}\right)\\
&&+D_{u}\frac{\triangle t}{\triangle y^{2}}\left(u_{i,j-1}^{n}-2u_{i,j}^{n}+u_{i,j+1}^{n}\right)\nonumber \\
&&-\triangle t\triangle_{i}^{n}\left(f^x_{i,j}\right)-\triangle t\triangle_{j}^{n}\left(f^y_{i,j}\right). \nonumber
\end{eqnarray}

The numerical method in the interior points of the 2D domain for the chemoattractant reads as:
\begin{eqnarray}\label{parabolic_chemo2D_discr_complete}
\qquad\qquad\phi_{i,j}^{n+1}&=&\phi_{i,j}^{n}+D_{\phi}\frac{\triangle t}{\triangle x^{2}}\left(\phi_{i-1,j}^{n}-2\phi_{i,j}^{n}+\phi_{i+1,j}^{n}\right)\\
&& +D_{\phi}\frac{\triangle t}{\triangle y^{2}}\left(\phi_{i,j-1}^{n}-2\phi_{i,j}^{n}+\phi_{i,j+1}^{n}\right)+\triangle t a u_{i,j}^{n}-\triangle t b\phi_{i,j}^{n} \nonumber.
\end{eqnarray}  
 
In the following we present the discretization of the boundary and transmission conditions to complete the numerical schemes. 
 
\subsubsection{Discretization of the boundary conditions for the doubly-parabolic problem}

Now, in order to complete our numerical scheme, we need to discretize the boundary conditions to obtain values for the boundary on each domain for the time step $n+1$.\\
   Since a qualitative characteristic of this model is the preservation of total mass, we want our numerical model to preserve mass at each time step. To this aim, we have to choose discrete boundary conditions that both are consistent with the analytical boundary conditions and preserve the mass in the numerical method. 
We remark that we present the computations without source term $g$ and we will add it in the sequel to complete the equations.\\

\vspace{0.2cm}

\underline{Boundary conditions for the density of individuals $u$.}\\
   The mass conservation over time on $\Omega_l$ reads as:
   \begin{eqnarray}
   I(t)=\int_{\Omega_l}u\left(t,x,y\right)d\Omega_l=\int_{\Omega_l}u\left(0,x,y\right)d\Omega_l=I(0).
   \end{eqnarray}
   Now, applying a quadrature rule for the numerical integration:
   \begin{eqnarray}
   \mathcal{I}^{n}\approx \int_{\Omega} u\left(t,x,y\right)d\Omega,
   \end{eqnarray}
   we need to ensure that 
    \begin{eqnarray}
   \mathcal{I}^{n+1}= \mathcal{I}^{n}.
   \end{eqnarray}
   For the numerical integration different quadrature formulas can be used.
   Since we want to use constant space-steps and want to obtain mass-preserving boundary conditions for the numerical methods, closed Newton-Cotes methods are suitable. In particular, we use the trapezoidal rule which introduces an integration error of $\mathcal{O}\left(\triangle x^{2}\right)$.\\
   For the one-dimensional trapezoidal rule with a function
   $z:\mathbb{R}\longrightarrow\mathbb{R}$ we have
   \begin{eqnarray}\label{trap1D}
   \int_{\Omega}z\left(x\right)d\Omega\approx \triangle x\left(\frac{F(x_{0})}{2}+\sum_{i=1}^{N}F(x_{i})+\frac{F(x_{N+1})}{2}\right)
   \end{eqnarray}
   and for the two-dimensional trapezoidal rule with function $z:\mathbb{R}^{2}\longrightarrow\mathbb{R}$ we have:
   \begin{eqnarray}\label{trap2D}
   \begin{array}{cl}
   \displaystyle\int_{\Omega}F(x,y)d\Omega_l&\approx 
    \frac{\triangle x\triangle y}{4}\Big(F(x_{0},y_{0})+F(x_{N_{x}+1},y_{0})+F(x_{0},y_{N_{y}+1})\\
    &+F(x_{N_{x}+1},y_{N_{y}+1})+2\displaystyle\sum_{i=1}^{N_{x}}\left(F(x_{i},y_{0})+F(x_{i},y_{N_{y}+1})\right)\\
   &+2\displaystyle\sum_{j=1}^{N_{y}}\left(F(x_{0},y_{j})+F(x_{N_{x}+1},y_{j})\right)+4\displaystyle\sum_{i=1}^{N_{x}}\displaystyle\sum_{j=1}^{N_{y}}F(x_{i},y_{j})\Big).
   \end{array}
   \end{eqnarray}
   
   Imposing the equality $\mathcal{I}^{n+1} - \mathcal{I}^{n} = 0$ in the 1D case gives:
   \begin{eqnarray}
    \triangle x\left(\frac{u_{0}^{n+1}}{2}-\frac{u_{0}^{n}}{2}+\displaystyle\sum_{i=1}^{N}\left(u_{i}^{n+1}-u_{i}^{n}\right)+\frac{u_{N+1}^{n+1}}{2}-\frac{u_{N+1}^{n}}{2}\right)=0.
    \end{eqnarray}
    Using the numerical scheme (\ref{uparabolic1d_discr}) for $u_{i}^{n+1}$ for $i=1,\dots,N$ we get:
    \begin{eqnarray*}
    \begin{array}{lllc}
  \triangle x&\left(\frac{u_{0}^{n+1}-u_{0}^{n}}{2} \right. &+D_{u_c}\frac{\triangle t}{\triangle x^{2}} \underbrace{\displaystyle\sum_{i=1}^{N}\left(u_{i-1}^{n}-2u_{i}^{n}+u_{i+1}^{n}\right)}_{=u_{0}^{n}-u_{1}^{n}-u_{N}^{n}+u_{N+1}^{n}} \\
  &&\left. -\triangle t\displaystyle\sum_{i=1}^{N}\Delta_{i}^{n}(f_{i})+\frac{u_{N+1}^{n+1}-u_{N+1}^{n}}{2} \right)=0.
  \end{array}
     \end{eqnarray*}
     
     Then we obtain:
      \begin{eqnarray*}
    \begin{array}{llcl}
   & u_{0}^{n+1}-u_{0}^{n} -2D_{u_c}\frac{\triangle t}{\triangle x^{2}}\left(u_{1}^{n}-u_{0}^{n}\right)
     +u_{N+1}^{n+1}-u_{N+1}^{n}
   \\&-2D_{u_c}\frac{\triangle t}{\triangle x^{2}}\left(u_{N}^{n}-u_{N+1}^{n}\right)-2\triangle t\displaystyle\sum_{i=1}^{N}\Delta_{i}^{n}(f_{i})=0.
     \end{array}
     \end{eqnarray*}
    
     Then, applying (\ref{star}) we obtain in the above formula:
     \begin{eqnarray}
     \begin{array}{lc}
   u_{0}^{n+1}-u_{0}^{n}-2D_{u_c}\frac{\triangle t}{\triangle x^{2}}\left(u_{1}^{n}-u_{0}^{n}\right)+\frac{\triangle t}{\triangle x}\left(f_{0}^{n}+f_{1}^{n}\right)\\
    +u_{N+1}^{n+1}-u_{N+1}^{n}-2D_{u_c}\frac{\triangle t}{\triangle x^{2}}\left(u_{N}^{n}-u_{N+1}^{n}\right)-\frac{\triangle t}{\triangle x}\left(f_{N}^{n}+f_{N+1}^{n}\right)=0. \nonumber
    \end{array}
     \end{eqnarray}
     We can now compute the values for both $u_{0}^{n+1}$ and $u_{N+1}^{n+1}$ so that the term equals to zero. By collecting values from nearby stencils together (otherwise we obtain an error of $\mathcal{O}(\triangle x)$ which can be verified by Taylor expansion), we obtain the following conditions at the outer boundaries of 1D domain: 
     \begin{equation}   \label{BC_parabolic1Dnumerical} 
     u_{0}^{n+1}=u_{0}^{n}+2D_{u_c}\frac{\triangle t}{\triangle x^{2}}\left(u_{1}^{n}-u_{0}^{n}\right)-\frac{\triangle t}{\triangle x}\left(f_{0}^{n}+f_{1}^{n}\right)
     \end{equation}
and     
     \begin{equation} \label{BC_parabolic1Dnumerical_dx} 
     u_{N+1}^{n+1} = u_{N+1}^{n}+2D_{u_c}\frac{\triangle t}{\triangle x^{2}}\left(u_{N}^{n}-u_{N+1}^{n}\right)+\frac{\triangle t}{\triangle x}\left(f_{N}^{n}+f_{N+1}^{n}\right).
     \end{equation}
       
     In the interior points of 1D domain we have the numerical scheme:
     \begin{equation}\label{parabolic1d_discr_complete}
     u_{i}^{n+1}=u_{i}^{n}+D_{u_c}\frac{\triangle t}{\triangle x^{2}}\left(u_{i-1}^{n}-2u_{i}^{n}+u_{i+1}^{n}\right)-\frac{\triangle t}{2\triangle x}\left(f_{i+1}^{n}-f_{i-1}^{n}\right),\qquad i=1,\ldots,N.
       \end{equation}
     Then we can state the following result.
\begin{Proposition}   
  The scheme (\ref{parabolic1d_discr_complete}) endowed with boundary conditions (\ref{BC_parabolic1Dnumerical}) and (\ref{BC_parabolic1Dnumerical_dx}) is mass-preserving by construction, since it is obtained imposing $\mathcal{I}^{n+1} - \mathcal{I}^{n} = 0$, as shown above.
  Moreover, the scheme, obtained with the integral method above, is second order in space up to the boundaries since it can be equivalently obtained using the following second-order approximation of the first derivative including a ghost cell:
     \begin{equation}\label{BCghostcell}
     \partial_{x}u(0)\approx \frac{u_{1}-u_{-1}}{2\triangle x}.
     \end{equation}
     
      Finally, the scheme is also positivity-preserving under the parabolic CFL condition.

      \end{Proposition}
     
     
     Now we compute $u_{0}^{n+1}$ directly by using the second-order centered numerical scheme and replace the ghost value $u_{-1}$ from the discretization of condition (\ref{BCghostcell}).
     While this works well when $f=0$, the same does not happen for $f\neq 0$, thus making the approach with the discrete integral equation still necessary.
     Futhermore, by using a different numerical integration scheme, we can achieve different mass-preserving boundary conditions of higher order.\\
Using the mass-preserving property argument, we compute boundary conditions for the corners and top and bottom boundaries of the 2D domain $\Omega_l$ for $f=0$.
By applying them with the numerical method (\ref{parabolic2d_discr_complete}) into $\mathcal{I}^{n+1}-\mathcal{I}^{n}=0$, we get the expression:
$$
\frac{\triangle t \triangle x}{4}\left(-4\triangle t\displaystyle\sum_{i=1}^{N_{x}}\displaystyle\sum_{j=1}^{N_{y}}\left(\Delta_{i}^{n}\left(f^x_{i,j}\right)+\Delta_{j}^{n}\left(f^y_{i,j}\right)\right)\right)=0,
$$
since the terms in $u$ cancel.
By choosing again the central in space second order finite difference scheme (\ref{star2D}) for $\text{div}\left(f\right)$, we get

\begin{eqnarray*}
\begin{array}{clc}
&\frac{1}{\triangle y}\displaystyle\sum_{i=1}^{N_{x}}\left(f_{i,N_{y}+1}^{y,n}+ f_{i,N_{y}}^{y,n}-f_{i,1}^{y,n}- f_{i,0}^{y,n}\right)\\
+&\frac{1}{\triangle x}\displaystyle\sum_{j=1}^{N_{y}}\left(f_{N_{x}+1,j}^{x,n}+ f_{N_{x},j}^{x,n}- f_{1,j}^{x,n}- f_{0,j}^{x,n}\right)=0.
\end{array}
\end{eqnarray*}

Now we can distribute the remaining values to the boundary values in the same way we did for the 1D-parabolic case. Therefore, we obtain the following mass-preserving boundary conditions:
\begin{equation}\label{BC_corners_2Dparabolic_source}
\left\{\begin{array}{llllll}
     u_{0,0}^{n+1}&=&u_{0,0}^{n}+2D_{u} \frac{\triangle t}{\triangle x^{2}}\left(u_{1,0}^{n}-u_{0,0}^{n}\right)+2D_{u}\frac{\triangle t}{\triangle y^{2}}\left(u_{0,1}^{n}-u_{0,0}^{n}\right)\\ 
     &&-\frac{\triangle t}{\triangle x}\left(f_{0,0}^{x,n}+f_{1,0}^{x,n}\right)-\frac{\triangle t}{\triangle y}\left(f_{0,0}^{y,n}+ f_{0,1}^{y,n}\right)\\&&+\triangle t g\left(x_{0},y_{0},t^{n},u_{0,0}^{n}\right)\\    
      u_{N_{x}+1,0}^{n+1}&=&u_{N_{x}+1,0}^{n}+2D_{u}\frac{\triangle t}{\triangle x^{2}}\left(u_{N_{x},0}^{n}-u_{N_{x}+1,0}^{n}\right)\\&&+2D \frac{\triangle t}{\triangle y^{2}}\left(u_{N_{x}+1,1}^{n}-u_{N_{x}+1,0}^{n}\right)\\
      &&-\frac{\triangle t}{\triangle x}\left(f_{N_{x}+1,0}^{x,n}+ f_{N_{x},0}^{x,n}\right)-\frac{\triangle t}{\triangle y}\left(f_{N_{x}+1,0}^{y,n}+ f_{N_{x}+1,1}^{y,n}\right)\\&&+\triangle t g\left(x_{N_{x}+1},y_{0},t_{n},u_{N_{x}+1,0}^{n}\right)\\
       u_{0,N_{y}+1}^{n+1}&=&u_{0,N_{y}+1}^{n} + 2D_{u} \frac{\triangle t}{\triangle x^{2}}\left(u_{1,N_{y}+1}^{n}-u_{0,N_{y}+1}^{n}\right)\\&&+2D_{u} \frac{\triangle t}{\triangle y^{2}}\left(u_{0,N_{y}}^{n}-u_{0,N_{y}+1}^{n}\right)\\
       &&-\frac{\triangle t}{\triangle x}\left(f_{0,N_{y}+1}^{x,n}+ f_{1,N_{y}+1}^{x,n}\right)-\frac{\triangle t}{\triangle y}\left(f_{0,N_{y}+1}^{y,n}+ f_{0,N_{y}}^{y,n}\right)\\&&+\triangle t g\left(x_{0},y_{N_{y}+1},t_{n},u_{0,N_{y}+1}^{n}\right)\\
        u_{N_{x}+1,N_{y}+1}^{n+1}&=&u_{N_{x}+1,N_{y}+1}^{n}+2D_{u} \frac{\triangle t}{\triangle x^{2}}\left(u_{N_{x},N_{y}+1}^{n}-u_{N_{x}+1,N_{y}+1}^{n}\right)
       \\ &&+2D_{u} \frac{\triangle t}{\triangle y^{2}}\left(u_{N_{x}+1,N_{y}}^{n}-u_{N_{x}+1,N_{y}+1}^{n}\right)\\
        &&-\frac{\triangle t}{\triangle x}\left(f_{N_{x}+1,N_{y}+1}^{x,n}+ f_{N_{x},N_{y}+1}^{x,n}\right)\\&&-\frac{\triangle t}{\triangle y}\left(f_{N_{x}+1,N_{y}+1}^{y,n}+ f_{N_{x}+1,N_{y}}^{y,n}\right)\\
\end{array}\right.
\end{equation}     
and for the top and bottom boundaries we have:
\begin{equation}\label{BC_borders_2Dparabolic_source}
\left\{\begin{array}{llllll}
u_{i,0}^{n+1}&=&u_{i,0}^{n}+D \frac{\triangle t}{\triangle x^{2}}\left(u_{i-1,0}^{n}-2u_{i,0}^{n}+u_{i+1,0}^{n}\right)+2D_{u} \frac{\triangle t}{\triangle y^{2}}\left(u_{i,1}^{n}-u_{i,0}^{n}\right)\\
&&-\frac{\triangle t}{2\triangle x}\left(f_{i+1,0}^{x,n}-f_{i-1,0}^{x,n}\right)-\frac{\triangle t}{\triangle y}\left(f_{i,0}^{y,n}+f_{i,1}^{y,n}\right)+\triangle t g\left(x_{i},y_{0},t_{n},u_{i,0}^{n}\right)\\
u_{i,N_{y}+1}^{n+1}&=&u_{i,N_{y}+1}^{n}+D \frac{\triangle t}{\triangle x^{2}}\left(u_{i-1,N_{y}+1}^{n}-2u_{i,N_{y}+1}^{n}+u_{i+1,N_{y}+1}\right)\\&&+2D \frac{\triangle t}{\triangle y^{2}}\left(u_{i,N_{y}}^{n}-u_{i,N_{y}+1}^{n}\right)\\
&&-\frac{\triangle t}{2\triangle x}\left(f_{i+1,N_{y}+1}^{x,n}- f_{i-1,N_{y}+1}^{x,n}\right)+\frac{\triangle t}{\triangle y}\left(f_{i,N_{y}}^{y,n}+ f_{i,N_{y}+1}^{y,n}\right)\\&&+\triangle t g\left(x_{i},y_{N_{y}+1},t_{n},u_{i,N_{y}+1}^{n}\right)\\
u_{0,j}^{n+1}&=&u_{0,j}^{n}+2D_{u} \frac{\triangle t}{\triangle x^{2}}\left(u_{1,j}^{n}-u_{0,j}^{n}\right)+D_{u}\frac{\triangle t}{\triangle y^{2}}\left(u_{0,j-1}^{n}-2u_{0,j}^{n}+u_{0,j+1}^{n}\right)\\
&&-\frac{\triangle t}{\triangle x}\left(f_{0,j}^{x,n}+ f_{1,j}^{x,n}\right)-\frac{\triangle t}{2\triangle y}\left(f_{0,j+1}^{y,n}- f_{0,j-1}^{y,n}\right)+\triangle t g\left(x_{0},y_{j},t_{n},u_{0,j}^{n}\right)\\
u_{N_{x}+1,j}^{n+1}&=&u_{N_{x}+1,j}^{n}+2D_{u} \frac{\triangle t}{\triangle x^{2}}\left(u_{N_{x},j}^{n}-u_{N_{x}+1,j}^{n}\right)\\&&+D_{u} \frac{\triangle t}{\triangle y^{2}}\left(u_{N_{x}+1,j-1}^{n}-2u_{N_{x}+1,j}^{n}+u_{N_{x}+1,j+1}^{n}\right)\\
&&-\frac{\triangle t}{2\triangle x}\left(f_{i+1,N_{y}+1}^{x,n}- f_{i.1,N_{y}+1}^{x,n}\right)+\frac{\triangle t}{\triangle x}\left(f_{i,N_{y}}^{y,n}+ f_{i,N_{y}+1}^{y,n}\right).\\
\end{array}\right.
\end{equation}     
     
     \vspace{0.8cm}
     \underline{Boundary conditions for the density of chemoattractant $\phi$.}\\
    
     For the computation of the conditions at the outer boudaries for the chemoattractant $\phi_{c}$ in the 1D-doubly parabolic model we proceed as above, but neglecting the source term $a_{c}u-b_{c}\phi_{c}$ to obtain boundary conditions that are mass-preserving.\\
     By doing so, we achieve the following second-order accurate and mass and positivity preserving boundary conditions for the chemoattractant:
     \begin{eqnarray}\label{BC_parabolic_chemo1D}
    \  & \phi_{0}^{n+1} = \phi_{0}^{n}+2D_{\phi_{c}}\frac{\triangle t}{\triangle x^{2}}\left(\phi_{1}^{n}-\phi_{0}^{n}\right)+\triangle t a_{c}u_{0}^{n}-\triangle tb_{c}\phi_{0}^{n}\\
     \ & \phi_{N+1}^{n+1}=\phi_{N+1}^{n}+2D_{\phi_{c}}\frac{\triangle t}{\triangle x^{2}}\left(\phi_{N}^{n}-\phi_{N+1}^{n}\right)+\triangle t a_{c}u_{N+1}^{n}-\triangle t b_{c}\phi_{N+1}^{n}.
     \end{eqnarray}
   The parabolic equation in the interior points is solved using a finite differences scheme in space and an explicit method in time:
    \begin{equation}\label{parabolic_chemo1D_discr}
     \phi_{i}^{n+1}=\phi_{i}^{n}+D_{\phi_{c}}\frac{\triangle t}{\triangle x^{2}}\left(\phi_{i-1}^{n}-2\phi^n_{i}+\phi_{i+1}^{n}\right)+\triangle t a_{c}u_{i}^{n}-\triangle t b_{c}\phi_{i}^{n}, \ i=1,\ldots,N.
     \end{equation}
     
In a similar way we can extend the numerical boundary conditions for the 2D-parabolic model.\\
   
     
Reasoning as above, we obtain the following boundary condition for the chemoattractant at the corners:
\begin{equation}\label{BC_corners2D_parabolic chemo}
\left\{\begin{array}{llllll}
\phi_{0,0}^{n+1}&=&\phi_{0,0}^{n}+2D_{\phi}\frac{\triangle t}{\triangle x^{2}}\left(\phi_{1,0}^{n}-\phi_{0,0}^{n}\right)+2D_{\phi}\frac{\triangle t}{\triangle y^{2}}\left(\phi_{0,1}^{n}-\phi_{0,0}^{n}\right)\\&&+\triangle t au_{0,0}^{n}-\triangle t b\phi_{0,0}^{n}\\
\phi_{N_{x}+1,0}^{n+1}&=&\phi_{N{x}+1,0}^{n}+2D_{\phi}\frac{\triangle t}{\triangle x^{2}}\left(\phi_{N_{x},0}^{n}-\phi_{N_{x}+1,0}^{n}\right)\\&&+2D_{\phi}\frac{\triangle t}{\triangle y^{2}}\left(\phi_{N_{x}+1,1}^{n}-\phi_{N_{x}+1,0}^{n}\right)+\triangle t a u_{N_{x}+1,0}^{n}-\triangle t b\phi_{N_{x}+1,0}^{n}\\
\phi_{N_{x}+1,N_{y}+1}^{n+1}&=&\phi_{N{x}+1,N_{y}+1}^{n}+2D_{\phi}\frac{\triangle t}{\triangle x^{2}}\left(\phi_{N_{x},N_{y}+1}^{n}-\phi_{N_{x}+1,N_{y}+1}^{n}\right)\\
&&+2D_{\phi}\frac{\triangle t}{\triangle y^{2}}\left(\phi_{N_{x}+1,N_{y}}^{n}-\phi_{N_{x}+1,N_{y}+1}^{n}\right)\\
&&+\triangle t a u_{N_{x}+1,N_{y}+1}^{n}-\triangle t b\phi_{N_{x}+1,N_{y}+1}^{n}\\
\phi_{0,N_{y}+1}^{n+1}&=&\phi_{0,0}^{n}+2D_{\phi}\frac{\triangle t}{\triangle x^{2}}\left(\phi_{1,N_{y}+1}^{n}-\phi_{0,N_{y}+1}^{n}\right)\\&&+2D_{\phi}\frac{\triangle t}{\triangle y^{2}}\left(\phi_{0,N_{y}}^{n}-\phi_{0,N_{y}+1}^{n}\right)+\triangle t a u_{0,N_{y}+1}^{n}-\triangle t b\phi_{0,N_{y}+1}^{n},
\end{array}\right.
\end{equation}
and for the borders we have:
\begin{equation}\label{BC_borders2D_parabolic_chemo}
\left\{\begin{array}{llllll}
\phi_{i,0}^{n+1}&=&\phi_{i,0}^{n}+D_{\phi}\frac{\triangle t}{\triangle x^{2}}\left(\phi^n_{i-1,0}-2\phi^n_{i,0}+\phi_{i+1,0}^{n}\right)+2D_{\phi}\frac{\triangle t}{\triangle y^{2}}\left(\phi_{i,1}^{n}-\phi_{i,0}^{n}\right)\\&&+\triangle t au_{i,0}^{n}-\triangle tb\phi_{i,0}^{n}\\
\phi_{i,N_{y}+1}^{n+1}&=&\phi_{i,N_{y}+1}^{n}+D_{\phi}\frac{\triangle t}{\triangle x^{2}}\left(\phi_{i-1,N_{y}+1}-2\phi_{i,N_{y}+1}+\phi_{i+1,N_{y}+1}^{n}\right)\\
&&+2D_{\phi}\frac{\triangle t}{\triangle y^{2}}\left(\phi_{i,N_{y}}^{n}-\phi_{i,N_{y}+1}^{n}\right)+\triangle t au_{i,N_{y}+1}^{n}-\triangle t b\phi_{i,N_{y}+1}^{n}\\
\phi_{0,j}^{n+1}&=&\phi_{0,j}^{n}+2D_{\phi}\frac{\triangle t}{\triangle x^{2}}\left(\phi^n_{1,j}-\phi^n_{0,j}\right)+D_{\phi}\frac{\triangle t}{\triangle y^{2}}\left(\phi^n_{0,j-1} - \phi_{0,j}^{n}+ \phi_{0,j+1}^{n}\right) \\
&&+\triangle t au_{0,j}^{n}-\triangle t b\phi_{0,j}^{n},\\
\phi_{N_{x}+1,j}^{n+1}&=&\phi_{N_{x}+1,j}^{n}+2D_{\phi}\frac{\triangle t}{\triangle x^{2}}\left(u_{N_{x},j}-u_{N_{x}+1,j}^{n}\right)\\
&&+D_{\phi}\frac{\triangle t}{\triangle y^{2}}\left(\phi_{N_{x}+1,j-1}^{n}-\phi_{N_{x}+1,j}^{n}+ \phi_{0,j+1}^{n}\right)\\
&&+\triangle t au_{N_{x}+1,j}^{n}-\triangle t b\phi_{N_{x}+1,j}^{n}.
\end{array}\right.
\end{equation}

We have now have a complete numerical method to solve (\ref{parabolic1d}) and (\ref{parabolic2d}).

\subsubsection{Discretization of the transmission conditions for the doubly-parabolic case}\label{sez:discrtrans}

Since for $K=0$ we would achieve the same analytical boundary conditions for separate domains (\ref{BCcells_analytical}) and (\ref{BCchemo_analytical}), we follow the same approach by using ghost values as in (\ref{BCghostcell}) in the numerical scheme to obtain the boundary conditions, since in such a way mass preserving and positivity preserving boundary condition are achieved.\\
By using the approximation formula (\ref{star}) in the condition (\ref{BCKK1_2D1Dparabolic}) on the left of node $1L$ we have:
\begin{eqnarray*}
\begin{array}{lccl}\label{BCKK2_2D1Dparabolic_ghostvalue}
&D_{u}\partial_{x} u \left(L_{x},y,t\right)-f^x \left(L_{x},y,t\right)&=&K\left(u_{c}\left(0,t\right)-u\left(L_{x},y,t\right)\right)\textrm{ for  } y\in\left[a_{1},b_{1}\right].
\end{array}
\end{eqnarray*}
Then we have:
\begin{eqnarray*}
\begin{array}{lccl}
& D_{u}\frac{u_{N_{x}+2,j}^{n}-u_{N_{x},j}^{n}}{2\triangle x}&=&K\left(u_{0}^{n}-u_{N_{x}+1,j}^{n}\right)+ f_{N_{x}+1,j}^{x,n}
\end{array}
\end{eqnarray*}
and we get:
\begin{eqnarray}\label{ghostend}
\begin{array}{lccl}
u_{N_{x}+2,j}^{n}&=&u_{N_{x},j}^{n}+K\frac{2\triangle x}{D}\left(u_{0}^{n}-u_{N_{x}+1,j}^{n}\right)+\frac{2\triangle x}{D_{u}} f_{N_{x}+1,j}^{x,n}
\end{array}
\end{eqnarray}
for $j=j_{a1},\dots,j_{b1}$.\\

Moreover, using (\ref{star}) in (\ref{BCKK2_2D1Dparabolic}), we can write:
\begin{equation}\label{KK_ghostvalue_1Dparabolic}
D_{u_c}\frac{u_{1}^{n}-u_{-1}^{n}}{2\triangle x} = \displaystyle\int_{a_{1}}^{b_{1}}K\left(u_{0}^{n}-u\left(L_{x},y,t\right)\right)dy+f_{0}^{n}
\end{equation}
and then we obtain:
\begin{eqnarray*}
\begin{array}{lccl}
u_{-1}^{n}&=&u_{1}^{n}-\frac{2\triangle x}{D_{u_c}}\displaystyle\int_{a_{1}}^{b_{1}}K\left(u_{c}\left(0,t\right)-u \left(L_{x},y,t\right)\right)dy-\frac{2\triangle x}{D_{u_c}}f_{0}^{n},
\end{array}
\end{eqnarray*}
and we finally get the formula:
\begin{equation}\label{bc1dparab}
u_{-1}^{n}=u_{1}^{n}-K\frac{2\triangle x}{D_{u_c}}\sigma u_{0}^{n}+\frac{2\triangle x}{D_{u_c}}\displaystyle\int_{a_{1}}^{b_{1}}K u\left(L_{x},y,t\right)dy-\frac{2\triangle x}{D_{u_c}}f_{0}^{n}.
\end{equation}

We now use the Ansatz to apply the ghost values into the numerical scheme without specific chemotactic approximation (\ref{parabolic1d_discr_complete}) and (\ref{parabolic2d_discr_complete}), and use the discrete integral equation to determine the chemotactic term discretization.\\
Because we now not only need to conserve the mass in each domain, but in both connected ones, the must use the expanded discrete integral equation to compute the total mass over both domains.\\

Plugging the ghost values (\ref{ghostend}) and (\ref{bc1dparab}), respectively, into the numerical schemes (\ref{parabolic2d_discr_complete}) and (\ref{parabolic1d_discr_complete}), we get the conditions at the interface (node $1L$):
\begin{equation}\label{parab_trasm_2d}
\begin{array}{lcll}
u_{N_{x}+1,j}^{n+1}&=&u_{N_{x}+1,j}^{n}+2D_{u} \frac{\triangle t}{\triangle x^{2}}\left(u_{N_{x},j}^{n}-u_{N_{x}+1,j}^{n}\right)+2K\frac{\triangle t}{\triangle x}\left(u_{0}^{n}-u_{N_{x}+1,j}^{n}\right)\\
&&+D_{u} \frac{\triangle t}{\triangle y^{2}}\left(u_{N_{x}+1,j-1}^{n}-2u_{N_{x}+1,j}^{n}+u_{N_{x}+1,j+1}^{n}\right)\\
&&+2\frac{\triangle t}{\triangle x}f_{N_{x}+1,j}^{x,n} -\triangle t \Delta_{N_{x}+1}^{n}\left(f^x_{j}\right)-\triangle t\Delta_{j}^{n}\left(f^y_{N_{x}+1}\right),
 \end{array}
\end{equation}
and
\begin{equation}\label{parab_trasm_1d}
\begin{array}{lcll}
 u_{0}^{n+1} &=& u_{0}^{n}+2D_{u_c}\frac{\triangle t}{\triangle x^{2}}\left(u_{1}^{n}-u_{0}^{n}\right)-2K\frac{\triangle t}{\triangle x}\sigma u_{0}^{n}+2K\frac{\triangle t}{\triangle x}\displaystyle\int_{a_{1}}^{b_{1}}u\left(L_{x},y,t_{n}\right)dy\\
&&-2\frac{\triangle t}{\triangle x}f_{0}^{n}-\triangle t\Delta_{0}^{n}\left(f\right).
 \end{array}
\end{equation}

In particular, the conservation of the discrete total mass reads as:
\begin{equation}\label{trap1D2D}
\mathcal{I}_{\text{1D}}^{n+1}+\mathcal{I}_{\text{2D}}^{n+1}-\mathcal{I}_{\text{1D}}^{n}-\mathcal{I}_{\text{2D}}^{n}=0,
\end{equation}
and now, applying the conditions (\ref{parab_trasm_2d}) and (\ref{parab_trasm_1d}) with the other boundary conditions (\ref{BC_borders_2Dparabolic_source}) and (\ref{BC_parabolic1Dnumerical}) we get:  
\begin{eqnarray*}
\begin{array}{lllll}
&\triangle x \left(\right.-K\frac{\triangle t}{\triangle x}\sigma u_{0}^{n}+K\frac{\triangle t}{\triangle x}\displaystyle\int_{a_{1}}^{b_{1}}u\left(L_{x},y,t_{n}\right)dy-\frac{\triangle t}{\triangle x}f_{0}^{n}\\
&\left.-\frac{\triangle t}{2}\Delta_{0}^{n}\left(f_{c}\right) -\frac{\triangle t}{2\triangle x}\left(-f_{0}^{n}-f_{1}^{n}\right)\right)\\
+& \frac{\triangle x\triangle y}{4}\left(\right.\left.2\displaystyle\sum_{j=j_{a_{1}}}^{j_{b_{1}}}\left(2K\frac{\triangle t}{\triangle x}\left(u_{0}^{n}-u_{N_{x}+1,j}^{n}\right)+2\frac{\triangle t}{\triangle x}f^x_{N_{x}+1,j}-\triangle t\Delta_{N_{x}+1}^{n}\left(f^x_{j}\right)\right.\right.\\
&\left.\left.-\triangle t\Delta_{j}^{n}\left(f^y_{N_{x}+1}\right)\right)\right.\\
&-\frac{2\triangle t}{\triangle x}\displaystyle\sum_{j=j_{a_{1}}}^{j_{b_{1}}}\left(f_{N_{x}+1,j}^{x,n}+f_{N_{x},j}^{x,n}\right)-\frac{2\triangle t}{\triangle y}\displaystyle\sum_{j=j_{a_{1}}}^{j_{b_{1}}}\left(f_{N_{x}+1,j}^{y,n}+ f_{N_{x},j}^{y,n}\right)\Big)=0,
\end{array}
\end{eqnarray*}
and obtain the following transmission conditions 

\begin{equation}\label{BCparabolic1DKK}
\begin{array}{lcl}
u_{0}^{n+1}&=& \underbrace{u_{0}^{n}+2D_{u_c}\frac{\triangle t}{\triangle x^{2}}\left(u_{1}^{n}-u_{0}^{n}\right)-\frac{\triangle t}{\triangle x}\left(f_{0}^{n}+f_{1}^{n}\right)}_{\text{same as for BC without transmission condition}}\\
&&-2K\frac{\triangle t}{\triangle x}\sigma u_{0}^{n}+2K\frac{\triangle t}{\triangle x}\displaystyle\int_{a_{1}}^{b_{1}}u\left(L_{x},y,t_{n}\right)dy
\end{array}
\end{equation}
\begin{equation}\label{BCtransmission2Dparabolicsourceterm}
\begin{array}{lcl}
u_{N_{x}+1,j}^{n+1} &=& u_{N_{x}+1,j}^{n}+2D_{u}\frac{\triangle t}{\triangle x^{2}}\left(u_{N_{x},j}^{n}-u_{N_{x}+1,j}^{n}\right)\\&&+D_{u_c}\frac{\triangle t}{\triangle y^{2}}\left(u_{N_{x}+1,j-1}^{n}-2u_{N_{x}+1,j}^{n}+u_{N_{x}+1,j+1}^{n}\right)\\
&&+\frac{\triangle t}{\triangle x}\left(f_{N_{x}+1,j}^{x,n}+ f_{N_{x},j}^{x,n}\right)-\frac{\triangle t}{2\triangle y}\left(f_{N_{x}+1,j+1}^{y,n}- f_{N_{x}+1,j-1}^{y,n}\right)\\
&& +\underbrace{2K\frac{\triangle t}{\triangle x}\left(u_{0}^{n}-u_{N_{x}+1,j}^{n}\right)}_{\text{additional term for transmission condition,}} 
\end{array}
\end{equation}
for $j=j_{a_{1}},\dots,j_{b_{1}}$.\\
The integral expression of the density $u_{0}^{n+1}$ in (\ref{BCparabolic1DKK}) can be expressed with a numerical quadrature form, such as the trapezoidal rule, as in (\ref{trap1D}).
\\
Proceeding analogously as above, this approach leads to mass-preserving and positivity-preserving  transmission conditions for the chemoattractant $\phi$ as well. In particular, we have at the first and last endpoint, respectively:
\begin{eqnarray}\label{BC1DchemoKK}
\phi_{0}^{n+1}&=&\phi_{0}^{n}+2D_{\phi_{c}}\frac{\triangle t}{\triangle x^{2}}\left(\phi_{1}^{n}-\phi_{0}^{n}\right)+\triangle t a_{c}u_{0}^{n}-\triangle t b_{c}\phi_{0}^{n}\\&&-2K\frac{\triangle t}{\triangle x}\sigma \phi_{0}^{n}+2K\frac{\triangle t}{\triangle x}\displaystyle\int_{a_{1}}^{b_{1}}\phi\left(L_{x},y,t_{n}\right)dy \nonumber
\end{eqnarray}
and
\begin{eqnarray}\label{BC2DchemoKK}
\begin{array}{lcl}
\phi_{N_{x}+1,j}^{n+1}&=&\phi_{N_{x}+1,j}^{n}+2D_{\phi}\frac{\triangle t}{\triangle x^{2}}\left(\phi_{N_{x},j}^{n}-\phi_{N_{x}+1,j}^{n}\right)\\
&&+D_{\phi}\frac{\triangle t}{\triangle y^{2}}\left(\phi_{N_{x}+1,j-1}^{n}-2\phi_{N_{x}+1,j}^{n}+\phi_{N_{x}+1,j+1}^{n}\right)\\
&&+\triangle au_{N_{x}+1,j}^{n}-\triangle t b\phi_{N_{x}+1,j}^{n}+2K\frac{\triangle t}{\triangle x}\left(\phi_{0}^{n}-\phi_{N_{x}+1,j}^{n}\right).
\end{array}
\end{eqnarray}
\\
We have finally developed a complete numerical scheme to treat doubly-parabolic partial differential equations systems in two domains, 1D and 2D, connected through a node, which ensures the mass conservation and the positivity as the original PDE.\\


\subsection{The hyperbolic-parabolic case}
The second order AHO scheme on a line was introduced in \cite{NR} for the 1D hyperbolic system (\ref{hyper-gen}). Here, considering the presence of the source term $g$ on the right hand side of the equation for the density of cells, the AHO scheme reads as:
\begin{equation}
\left\{\begin{array}{lcl}\label{AHO2_scheme_cells}
u_{i}^{n+1}&=&u_{i}^{n}+\lambda\frac{\triangle t}{2\triangle x}\left(u_{i-1}^{n}-2u_{i}^{n}+u_{i+1}^{n}\right)-\left(\frac{\triangle t}{2\triangle x}-\frac{\triangle t}{4\lambda}\right)\left(v_{i+1}^{n}-v_{i-1}^{n}\right)\\
&&+\frac{\triangle t}{4\lambda}\left(f_{i-1}^{n}-f_{i+1}^{n}\right)+\\&&\frac{\triangle t}{4}\left(g\left(x_{i-1},t_{n},u_{i-1}^{n}\right)+2g\left(x_{i},t_{n},u_{i}^{n}\right)+g\left(x_{i+1},t_{n},u_{i+1}^{n}\right)\right),\\
\\
v_{i}^{n+1}&=&v_{i}^{n}-\lambda^{2}\frac{\triangle t}{2\triangle x}\left(u_{i+1}^{n}-u_{i-1}^{n}\right)+\frac{\lambda\triangle t}{2\triangle x}\left(v_{i-1}^{n}-2v_{i}^{n}+v_{i+1}^{n}\right)\\
&&-\frac{\triangle t}{4}\left(v_{i-1}^{n}+2v_{i}^{n}+v_{i+1}^{n}\right)+\frac{\triangle t}{4}\left(f_{i-1}^{n}+2f_{i}^{n}+f_{i+1}^{n}\right)\\&&+\lambda\frac{\triangle t}{4}\left(g\left(x_{i-1},t_{n},u_{i-1}^{n}\right)-g\left(x_{i+1},t_{n},u_{i+1}^{n}\right)\right),
\end{array}\right.
\end{equation}
with mass-preserving boundary conditions (including the additional source term $g$) at the external boundaries. 
We remark that for the hyperbolic-parabolic model not only  mass must be preserved as the in the fully-parabolic model, but also the flux $v$ needs to converge towards the steady state $v=0$. Since here we have the 1D domain connected at both the endpoints we do not need to use numerical boundary conditions for the outer boundaries. However, for the details and the description of the AHO numerical scheme at the outer boundaries, see \cite{NR}.\\
For this reason we use the so called AHO (Asymptotic Higher Order) schemes (see \cite{BNR14} for the study of AHO scheme at interfaces including mass-preserving transmission conditions) with source term $g$ for which the approximation of the stationary solutions is up to third order of accuracy and it converges towards a numerical solution with $v=0$, while preserving the mass.\\

\subsubsection{Discretization of transmission conditions for the 2D-doubly-parabolic and 1D-hyperbolic-parabolic case}\label{sec:trasmhyp}

The first equation is the same as for the 2D-doubly-parabolic and 1D-doubly-parabolic case. Hence we derive the same transmission condition for $u_{N_{x}+1,j}^{n+1}$ for $j=j_{a_{1}},\dots,j_{b_{1}}$.\\
For the flux, the transmission condition (\ref{BCKK2_2D1D_hyperbolic}) gives us
\begin{equation}\label{BCKK2_2D1D_hyperbolic_flux}
v_{0}^{n+1} = -K\sigma u_{0}^{n+1}+K\displaystyle\int_{a_{1}}^{b_{1}}u\left(L_{x},y,t_{n+1}\right)dy.
\end{equation}
We remark that here we have two problems  which do not occur in the previous model. In particular:
we do not have a transmission condition formula for $u_{0}^{n+1}$ and the formula for $v_{0}^{n+1}$ is implicit, since it depends on values at time step $t_{n+1}$.\\
In order to solve these issues, we will use once again the equivalence between discrete total masses to obtain mass-preserving computation formula for $u_{0}^{n+1}$ and, since $v_{0}^{n+1}$ only depends on $u_{0}^{n+1}$ and $u_{N_{x}+1,j}^{n+1}$, which can be obtained explicity from computed values at time step $t_n$, we can calculate $v_{0}^{n+1}$ explicity as well.\\
Then, imposing that:
$$
\mathcal{I}_{\text{2D}}^{n+1}+ \mathcal{I}_{\text{2D}}^{n+1}- \mathcal{I}_{\text{1D}}^{n}- \mathcal{I}_{\text{1D}}^{n}=0
$$
we get:
\begin{eqnarray*}
\begin{array}{llll}
&\frac{\triangle x}{2}\left[ u_{0}^{n+1}-u_{0}^{n}+\lambda\frac{\triangle t}{\triangle x}\left(u_{0}^{n}-u_{1}^{n}\right)-\left(\frac{\triangle t}{\triangle x}-\frac{\triangle t}{2\lambda}\right)\left(-v_{0}^{n}-v_{1}^{n}\right)+\frac{\triangle t}{2\lambda}\left(f_{0}^{n}+f_{1}^{n}\right)\right]&&\\
&+\frac{\triangle x\triangle y}{4}\left[4\displaystyle\sum_{j=j_{a_{1}}}^{j_{b_{1}}}\frac{\triangle t K}{\triangle x}\left(u_{0}^{n}-u_{N_{x}+1,j}^{n}\right)\right] = 0
\end{array}
\end{eqnarray*}
and we finally obtain the  transmission condition with source term $g$ :
\begin{eqnarray}\label{BCtransmission_1Dhyperbolic}
\begin{array}{lcl}
\Longrightarrow u_{0}^{n+1}&=&u_{0}^{n}+\lambda\frac{\triangle t}{\triangle x}\left(u_{1}^{n}-u_{0}^{n}\right)-\left(\frac{\triangle t}{\triangle x}-\frac{\triangle t}{2\lambda}\right)\left(v_{0}^{n}+v_{1}^{n}\right)-\frac{\triangle t}{2\lambda}\left(f_{0}^{n}+f_{1}^{n}\right)\\
&&+\frac{\triangle t}{2}\left(g\left(x_{0},t_{n},u_{0}^{n}\right)+g\left(x_{1},t_{n},u_{1}\right)\right)\\&&-2K\frac{\triangle t}{\triangle x}\triangle y \displaystyle\sum_{j=j_{a_{1}}}^{j_{b_{1}}}\left(u_{0}^{n}-u_{N_{x}+1,j}^{n}\right).
\end{array}
\end{eqnarray}

\begin{Proposition}
The complete numerical scheme derived for the 2D-doubly-parabolic-1D-hyperbolic-parabolic model has the feature to be mass-preserving across the transmission conditions. Note that for the chemoattractant equation is the same as for the 1D-doubly-parabolic and 2D-doubly-parabolic case. Hence, the numerical schemes (\ref{parabolic_chemo1D_discr}) and (\ref{parabolic_chemo2D_discr_complete}) with boundary conditions (\ref{BC_parabolic_chemo1D}) and (\ref{BC2DchemoKK}) can be used.
\end{Proposition}

\subsubsection{Multiple channels.}\label{multi}
In the previous paragraphs we have connected the two-dimensional domain $\Omega_l$ with a single one-dimensional channel $I$ at $(L_{x},y)\in\Omega_l$ with $ y \in\left[a_{1},b_{1}\right]$,  and  $j_{a_{1}}$ and $j_{b_{1}}$, the positions of the endpoints of the corridor on the numerical grid. Of course this can be extended to more channels.\\
Let $\left(I_m\right)_{m=1,\dots,M}$ be $\mathcal{M}$ corridors, connected to the two-dimensional domain $\Omega_l$ at $\left(L_{x},y_{m}\right)$ with $y_{m}\in\left[a_{m},b_{m}\right]$ and $a_1>0$, $b_{m}<a_{m+1}$, for $m=1,\dots,\mathcal{M}-1$, and $b_{\mathcal{M}}<L_{y}$ to avoid intersections of the corridors, with equal width  $\sigma:=b_{m}-a_{m}=k \triangle y$, $k\in\mathbb{N}$.\\
\subsubsection{Implemented algorithm.}\label{sec:algo}
Before presenting the numerical tests in the next section \ref{sec:tests}, we detail the approximation scheme for the density $u$, also including the source term $g$, implemented to solve the problem in the 2D-1D domain.
As underlined before, it is necessary to use implicit schemes to consider the presence of stiff source terms. For this reason, for the approximation of the time derivatives we use the Crank-Nicolson method on the diffusion and source term, which is a second order implicit method and the explicit central method for the convection term.\\
Because of the explicit term, we have numerical restrictions on the mesh grid and time step. Furthermore, as discussed previously, we introduce artificial viscosity to avoid oscillations due to not suitable mesh grid size in dominant convection regime, which is often the case in chemotaxis models.
The implicit-explicit numerical method used to compute the solutions for the density $u$ in (\ref{parabolic2d}) inside the 2D domain $\Omega_l$ is:
\begin{equation}\label{CNmethod} 
\begin{array}{lcl}
u_{i,j}^{n+1}&=& u_{i,j}^{n} +D_{u}\frac{\triangle t}{2}\left[\frac{\delta_{x}^{2}(u_{i,j}^{n}+u_{i,j}^{n+1})}{\triangle x^{2}}+\frac{\delta_{y}^{2}(u_{i,j}^{n}+u_{i,j}^{n+1})}{\triangle y^{2}}\right] \\
&& -\frac{\triangle t}{4}  \biggl[\frac{\delta_{x}^{0}f_{i,j}^{n}}{\triangle x}+\frac{\delta_{y}^{0}f_{i,j}^{n}}{\triangle y}\biggr]+\frac{\triangle t}{2}\left( g_{i,j}^{n}+g_{i,j}^{n+1}\right)\\
 &&\underbrace{-\triangle t\biggl[\frac{\delta_{x}^{2}\theta_{i,j}^{n}}{2\triangle x}+\frac{\delta_{y}^{2}\theta_{i,j}^{n}}{2\triangle y}\biggr]}_{\text{artificial viscosity}},
\end{array}
\end{equation}
with $\theta_{i,j}^{n}:=\chi(u_{i,j}^{n},\varphi_{i,k}^{n})\vert\nabla\varphi_{i,j}^{n}\vert$ for $i=1,\dots,N_{x}, \ j=1,\dots,N_{y}$.
As can be seen, the function $\theta$ used for the artificial viscosity is almost identical to $f$ with the exception of using the absolute value of $\nabla \varphi$.
By using this, we increase artificial viscosity only where the gradient of the chemoattractant increases. This reduces the restriction on the meshgrid due to the condition induced by the cell P\'eclet number.
 The numerical transmission condition on the left of node $1L$ ($i=N_x+1, j=j_{a_1},\dots,j_{b_1}$) is:
\begin{equation}\label{BCtransmission2Dparabolicsourcetermg}
\begin{array}{lcl}
u_{N_{x}+1,j}^{n+1} &=& u_{N_{x}+1,j}^{n}-D_{u}\frac{\triangle t}{\triangle x^{2}}\delta_{x}^{1}(u_{N_{x},j}^{n}+u_{N_{x},j}^{n+1})+D_{u}\frac{\triangle t}{2\triangle y^{2}}\delta_{y}^{2}\left(u_{N_{x}+1,j}^{n}+u_{N_{x}+1,j}^{n+1}\right)\\
&&+\frac{\triangle t}{\triangle x}\left(f_{N_{x}+1,j}^{x,n}+ f_{N_{x},j}^{x,n}\right)-\frac{\triangle t}{\triangle y}\left(f_{N_{x}+1,j+1}^{y,n}- f_{N_{x}+1,j-1}^{y,n}\right)\\
&&+\frac{\triangle t}{2}\left(g(x_{N_x+1},y_j,t_n, u^n_{N_x+1,j})+g(x_{N_x+1},y_j,t_{n+1}, u^{n+1}_{N_x+1,j})\right)\\
&&-\triangle t\left(\frac{\delta_{x}^{1}\theta_{N_{x},j}^{n}}{\triangle x}+\frac{\delta_{y}^{2}\theta_{N_{x+1},j}^{n}}{2\triangle y}\right)\\
&& +\underbrace{K\frac{\triangle t}{\triangle x}\left(u_{0}^{n}-u_{N_{x}+1,j}^{n}+u_{0}^{n+1}-u_{N_{x}+1,j}^{n+1}\right)}_{\text{additional term for transmission condition}}. 
\end{array}
\end{equation}
The role of KK coefficient $K$ in the positivity of (\ref{BCtransmission2Dparabolicsourcetermg}) is discussed in Remark \ref{remark_kk}.
\newpage
For the corners we use the following boundary conditions:
\begin{equation}\label{BC_corners_2Dparabolic_sourceg}
\left\{\begin{array}{llllll}
     u_{0,0}^{n+1}&=&u_{0,0}^{n}+D_{u} \frac{\triangle t}{\triangle x^{2}}\delta_{x}^{1}\left(u_{0,0}^{n}+u_{0,0}^{n+1}\right)+D_{u}\frac{\triangle t}{\triangle y^{2}}\delta_{y}^{1}\left(u_{0,0}^{n}+u_{0,0}^{n+1}\right)\\ 
     &&-\frac{\triangle t}{\triangle x}\left(f_{0,0}^{x,n}+f_{1,0}^{x,n}\right)-\frac{\triangle t}{\triangle y}\left(f_{0,0}^{y,n}+ f_{0,1}^{y,n}\right)\\&&+\frac{\triangle t}{2}\left(g\left(x_{0},y_{0},t_{n},u_{0,0}^{n}\right)+g\left(x_{0},y_{0},t_{n+1},u_{0,0}^{n+1}\right)\right)\\
     &&-\triangle t\left(\frac{\delta_{x}^{1}\theta_{0,0}^{n}}{\triangle x}+\frac{\delta_{y}^{0}\theta_{0,0}^{n}}{\triangle y}\right)\\    
      u_{N_{x}+1,0}^{n+1}&=&u_{N_{x}+1,0}^{n}-D_{u}\frac{\triangle t}{\triangle x^{2}}\delta_{x}^{1}\left(u_{N_{x},0}^{n}+u_{N_{x},0}^{n+1}\right)\\&&+D_{u} \frac{\triangle t}{\triangle y^{2}}\delta_{y}^{1}\left(u_{N_{x}+1,0}^{n}+u_{N_{x}+1,0}^{n+1}\right)\\
      &&-\frac{\triangle t}{\triangle x}\left(f_{N_{x}+1,0}^{x,n}+ f_{N_{x},0}^{x,n}\right)-\frac{\triangle t}{\triangle y}\left(f_{N_{x}+1,0}^{y,n}+ f_{N_{x}+1,1}^{y,n}\right)\\&&+\frac{\triangle t}{2} \left(g\left(x_{N_{x}+1},y_{0},t_{n},u_{N_{x}+1,0}^{n}\right)+g\left(x_{N_{x}+1},y_{0},t_{n+1},u_{N_{x}+1,0}^{n+1}\right)\right)\\
      &&-\triangle t\left(\frac{\delta_{x}^{1}\theta_{N_{x},0}^{n}}{\triangle x}+\frac{\delta_{y}^{1}\theta_{N_{x}+1,0}^{n}}{\triangle y}\right)\\
       u_{0,N_{y}+1}^{n+1}&=&u_{0,N_{y}+1}^{n} + D_{u} \frac{\triangle t}{\triangle x^{2}}\delta_{x}^{1}\left(u_{0,N_{y}+1}^{n}+u_{0,N_{y}+1}^{n+1}\right)\\&&-D_{u} \frac{\triangle t}{\triangle y^{2}}\delta_{y}^{1}\left(u_{0,N_{y}}^{n}+u_{0,N_{y}}^{n+1}\right)\\
       &&-\frac{\triangle t}{\triangle x}\left(f_{0,N_{y}+1}^{x,n}+ f_{1,N_{y}+1}^{x,n}\right) -\frac{\triangle t}{\triangle y}\left(f_{0,N_{y}+1}^{y,n}+ f_{0,N_{y}}^{y,n}\right)\\&&+\frac{\triangle t}{2}\left( g\left(x_{0},y_{N_{y}+1},t_{n},u_{0,N_{y}+1}^{n}\right)+g\left(x_{0},y_{N_{y}+1},t_{n+1},u_{0,N_{y}+1}^{n+1}\right)\right)\\
       &&-\triangle t\left(\frac{\delta_{x}^{1}\theta_{0,N_{y}+1}^{n}}{\triangle x}+\frac{\delta_{y}^{1}\theta_{0,N_{y}}^{n}}{\triangle y}\right)\\
        u_{N_{x}+1,N_{y}+1}^{n+1}&=&u_{N_{x}+1,N_{y}+1}^{n}-D_{u} \frac{\triangle t}{\triangle x^{2}}\delta_{x}^{1}\left(u_{N_{x},N_{y}+1}^{n}+u_{N_{x},N_{y}+1}^{n+1}\right)
       \\ &&-D_{u} \frac{\triangle t}{\triangle y^{2}}\delta_{y}^{1}\left(u_{N_{x}+1,N_{y}}^{n}+u_{N_{x}+1,N_{y}}^{n+1}\right)\\
        &&-\frac{\triangle t}{\triangle x}\left(f_{N_{x}+1,N_{y}+1}^{x,n}+ f_{N_{x},N_{y}+1}^{x,n}\right)\\&&-\frac{\triangle t}{\triangle y}\left(f_{N_{x}+1,N_{y}+1}^{y,n}+ f_{N_{x}+1,N_{y}}^{y,n}\right)\\
        &&+\frac{\triangle t}{2}\left( g\left(x_{N_{x}+1},y_{N_{y}+1},t_{n},u_{N_{x}+1,N_{y}+1}^{n}\right)\right.\\&&\left.+g\left(x_{N_{x}+1},y_{N_{y}+1},t_{n+1},u_{N_{x}+1,N_{y}+1}^{n+1}\right)\right)\\
        &&-\triangle t\left(\frac{\delta_{x}^{1}\theta_{N_{x},N_{y}+1}^{n}}{\triangle x}+\frac{\delta_{y}^{1}\theta_{N_{x}+1,N_{y}}^{n}}{\triangle y}\right)
\end{array}\right.
\end{equation}     
and for the top and bottom boundaries:
\begin{equation}\label{BC_borders_2Dparabolic_sourceg}
\left\{\begin{array}{llllll}
u_{i,0}^{n+1}&=&u_{i,0}^{n}+D_{u} \frac{\triangle t}{2\triangle x^{2}}\delta_{x}^{2}\left(u_{i,0}^{n}+u_{i,0}^{n+1}\right)+D_{u} \frac{\triangle t}{\triangle y^{2}}\delta_{y}^{1}\left(u_{i,0}^{n}+u_{i,0}^{n+1}\right)\\
&&-\frac{\triangle t}{2\triangle x}\delta_{x}^{0}\left(f_{i,0}^{x,n}\right)-\frac{\triangle t}{\triangle y}\left(f_{i,0}^{y,n}+f_{i,1}^{y,n}\right)\\&&+\frac{\triangle t}{2}\left( g\left(x_{i},y_{0},t_{n},u_{i,0}^{n}\right)+g\left(x_{i},y_{0},t_{n+1},u_{i,0}^{n+1}\right)\right)\\
&& -\triangle t\left(\frac{\delta_{x}^{2}\theta_{i,0}^{n}}{2\triangle x}+\frac{\delta_{y}^{1}\theta_{i,0}^{n}}{\triangle y}\right)\\
u_{i,N_{y}+1}^{n+1}&=&u_{i,N_{y}+1}^{n}+D_{u} \frac{\triangle t}{2\triangle x^{2}}\delta_{x}^{2}\left(u_{i,N_{y}+1}^{n}+u_{i,N_{y}+1}^{n+1}\right)\\&&-D_{u} \frac{\triangle t}{\triangle y^{2}}\delta_{y}^{1}\left(u_{i,N_{y}}^{n}+u_{i,N_{y}}^{n+1}\right)\\
&&-\frac{\triangle t}{2\triangle x}\delta_{x}^{0}\left(f_{i,N_{y}+1}^{x,n}\right)+\frac{\triangle t}{\triangle y}\left(f_{i,N_{y}}^{y,n}+ f_{i,N_{y}+1}^{y,n}\right)\\&&+\frac{\triangle t}{2}\left( g\left(x_{i},y_{N_{y}+1},t_{n},u_{i,N_{y}+1}^{n}\right)+ g\left(x_{i},y_{N_{y}+1},t_{n+1},u_{i,N_{y}+1}^{n+1}\right)\right)\\
&&-\triangle t\left(\frac{\delta_{x}^{2}\theta_{i,N_{y}+1}^{n}}{2\triangle x}+\frac{\delta_{y}^{1}\theta_{i,N_{y}}^{n}}{\triangle y}\right)\\
u_{0,j}^{n+1}&=&u_{0,j}^{n}+D_{u} \frac{\triangle t}{\triangle x^{2}}\delta_{x}^{1}\left(u_{0,j}^{n}+u_{0,j}^{n+1}\right)+D_{u}\frac{\triangle t}{2\triangle y^{2}}\delta_{y}^{2}\left(u_{0,j}^{n}+u_{0,j}^{n+1}\right)\\
&&-\frac{\triangle t}{\triangle x}\left(f_{0,j}^{x,n}+ f_{1,j}^{x,n}\right)
-\frac{\triangle t}{2\triangle y}\delta_{y}^{0}\left(f_{0,j}^{y,n}\right)\\&&+\frac{\triangle t}{2}\left(g\left(x_{0},y_{j},t_{n},u_{0,j}^{n}\right)+g\left(x_{0},y_{j},t_{n+1},u_{0,j}^{n+1}\right)\right)\\
&&-\triangle t\left(\frac{\delta_{x}^{1}\theta_{0,j}^{n}}{\triangle x}+\frac{\delta_{y}^{2}\theta_{0,j}^{n}}{2\triangle y}\right)\\
u_{N_{x}+1,j}^{n+1}&=&u_{N_{x}+1,j}^{n}-D_{u} \frac{\triangle t}{\triangle x^{2}}\delta_{x}^{1}\left(u_{N_{x},j}^{n}+u_{N_{x},j}^{n+1}\right)\\&&+D \frac{\triangle t}{2\triangle y^{2}}\delta_{y}^{2}\left(u_{N_{x}+1,j}^{n}+u_{N_{x}+1,j}^{n+1}\right)\\
&&-\frac{\triangle t}{2\triangle x}\delta_{x}^{0}\left(f_{i,N_{y}+1}^{x,n}\right)+\frac{\triangle t}{\triangle x}\left(f_{i,N_{y}}^{y,n}+ f_{i,N_{y}+1}^{y,n}\right)\\
&&+\frac{\triangle t}{2} \left(g\left(x_{N_{x}+1},y_{j},t_{n},u_{N_{x}+1,j}^{n}\right)+ g\left(x_{N_{x}+1},y_{j},t_{n}+1,u_{N_{x}+1,j}^{n+1}\right)\right)\\
&&-\triangle t\left(\frac{\delta_{x}^{1}\theta_{N_{x},j}^{n}}{\triangle x}+\frac{\delta_{y}^{2}\theta_{N_{x}+1,j}^{n}}{2\triangle y}\right).
\end{array}\right.
\end{equation}
Similarly, for the chemoattractant $\phi$ we have the implicit-explicit scheme in the interior points of the 2D domain:
\begin{equation}\label{CNmethod_phi}
\begin{array}{lcll}
\phi_{i,j}^{n+1}&=& \phi_{i,j}^{n} +D_{\phi}\frac{\triangle t}{2}\left[\frac{\delta_{x}^{2}\left(\phi_{i,j}^{n}+\phi_{i,j}^{n+1}\right)}{\triangle x^{2}}+\frac{\delta_{y}^{2}\left(\phi_{i,j}^{n}+ \phi_{i,j}^{n+1}\right)}{\triangle y^{2}}\right] \\
&& \frac{\triangle t}{2}(a (u^n_{i,j} + u^{n+1}_{i,j}) - \frac{\triangle t}{2}(b (\phi^n_{i,j} + \phi^{n+1}_{i,j}),
\end{array}
\end{equation}
and for the boundaries and the corners the numerical schemes for $\phi$ are, respectively, (\ref{BC_corners2D_parabolic chemo}) and (\ref{BC_borders2D_parabolic_chemo}).\\

\begin{remark}
If we consider two-dimensional domain $\Omega_{r}$ connected to the right endpoint of the one-dimensional corridor $I$, the complete numerical scheme for the left domain $\Omega_l$ described above can be considered.\\
The main difference is that the transmission conditions at the interface between the box and the channel (the left for the box $\Omega_l$ and the right for the corridor) are reversed to the left for the corridor and the right for the box $\Omega_{r}$. In the numerical scheme, the only change affects the channel $I$, where we have transmission conditions also for $u_{N+1}^{n}$ (resp. $v_{N+1}^{n}$). The same boundary condition can be used without transmission conditions, with only the additional term derived from the KK-condition and it must be added as well for $u_{0}^{n}$ (resp. $v_{0}^{n}$).
\end{remark} 
 For the computation of solutions on the one-dimensional channel $I$, we have two different approximations depending on the choice of the model we assign on it. If we solve the doubly-parabolic problem (\ref{parabolic1d}), the approximation scheme used is the Crank-Nicolson scheme, as above:
 \begin{equation}\label{CNmethod1d}
\begin{array}{lcll}
u_{i}^{n+1}&=& u_{i}^{n} +D_{u_c}\frac{\triangle t}{2}\left[\frac{\delta_{x}^{2}\left(u_{i}^{n}+u_{i}^{n+1}\right)}{\triangle x^{2}}\right]-\frac{\triangle t}{2}  \biggl[\frac{\delta_{x}^{0}\left(f_{i}^{n}\right)}{\triangle x} \biggr]+\frac{\triangle t}{2}\left( g_{i}^{n}+g_{i}^{n+1}\right)-\triangle t\left(\frac{\delta_{x}^{2}\theta_{i}^{n}}{2\triangle x}\right),\\
\end{array}
\end{equation}
 with the transmission condition on the left of node $1L$ ($i=0$) given by:
 \begin{equation}\label{BCparabolic1DKKsourceg}
 \begin{array}{lcl}
u_{0}^{n+1} &= &\underbrace{u_{0}^{n}+D_{u_c}\frac{\triangle t}{\triangle x^{2}}\delta_{x}^{1}\left(u_{0}^{n}+u_{0}^{n+1}\right)-\frac{\triangle t}{\triangle x}\left(f_{0}^{n}+f_{1}^{n}\right)-\triangle t\frac{\delta_{x}^{1}\theta_{0}^{n}}{\triangle x}}_{\text{same as for BC without transmission condition}}\\&& + \frac{\triangle t}{2}\left(g(x_0,t_n,u_0^n)+g(x_0,t_{n+1},u_0^{n+1})\right)\\
&&-K\frac{\triangle t}{\triangle x}\sigma (u_{0}^{n}+u_{0}^{n+1})+K\frac{\triangle t}{\triangle x}\displaystyle\sum_{j=j_{a_{1}}}^{j_{b_{1}}}\left(u_{N_{x}+1,j}^{n}+u_{N_{x}+1,j}^{n+1}\right).
\end{array}
\end{equation}

If, instead, we need to solve the hyperbolic-parabolic problem (\ref{hyperbolic1d}), an implicit version of the second order AHO scheme (\ref{AHO2_scheme_cells}) is used. Indeed, because of the different scales of the parameters, instabilities in the corridors can occur if the stability condition of the AHO scheme is not satisfied. In particular, inside the channels we use the scheme:
\begin{eqnarray*}
\left\{\begin{array}{lcl}\label{AHO2_scheme_cells_impl}
u_{i}^{n+1}&=&u_{i}^{n}+\lambda\frac{\triangle t}{2\triangle x}\left(u_{i-1}^{n+1}-2u_{i}^{n+1}+u_{i+1}^{n+1}\right)-\left(\frac{\triangle t}{2\triangle x}-\frac{\triangle t}{4\lambda}\right)\left(v_{i+1}^{n+1}-v_{i-1}^{n+1}\right)\\
&&+\frac{\triangle t}{4\lambda}\left(f_{i-1}^{n+1}-f_{i+1}^{n+1}\right)+\\&&\frac{\triangle t}{4}\left(g\left(x_{i-1},t_{n+1},u_{i-1}^{n+1}\right)+2g\left(x_{i},t_{n},u_{i}^{n+1}\right)+g\left(x_{i+1},t_{n+1},u_{i+1}^{n+1}\right)\right),\\
\\
v_{i}^{n+1}&=&v_{i}^{n}-\lambda^{2}\frac{\triangle t}{2\triangle x}\left(u_{i+1}^{n+1}-u_{i-1}^{n+1}\right)+\left(\frac{\lambda\triangle t}{2\triangle x}-\frac{\triangle t}{4}\right)\left(v_{i-1}^{n+1}-2v_{i}^{n+1}+v_{i+1}^{n+1}\right)\\
&&+\frac{\triangle t}{4}\left(f_{i-1}^{n+1}+2f_{i}^{n+1}+f_{i+1}^{n+1}\right)\\&&+\lambda\frac{\triangle t}{4}\left(g\left(x_{i-1},t_{n+1},u_{i-1}^{n+1}\right)-g\left(x_{i+1},t_{n+1},u_{i+1}^{n+1}\right)\right),
\end{array}\right.
\end{eqnarray*}
endowed with the following implicit version of transmission condition (\ref{BCKK2_2D1D_hyperbolic_flux}) :
\begin{eqnarray}\label{impl_transm_hyp}
u_{0}^{n+1}&=&u_{0}^{n}+\lambda\frac{\triangle t}{\triangle x}\left(u_{1}^{n+1}-u_{0}^{n+1}\right)-\left(\frac{\triangle t}{\triangle x}-\frac{\triangle t}{2\lambda}\right)\left(v_{0}^{n+1}+v_{1}^{n+1}\right) \nonumber\\
&&-\frac{\triangle t}{2\lambda}\left(f_{0}^{n+1}+f_{1}^{n+1}\right)+\frac{\triangle t}{2}\left(g\left(x_{0},t_{n+1},u_{0}^{n+1}\right)+g\left(x_{1},t_{n+1},u_{1}\right)\right)\\
&&-K\frac{\triangle t}{\triangle x}\triangle y \displaystyle\sum_{j=j_{a_{1}}}^{j_{b_{1}}}\left(u_{0}^{n+1}+u_{0}^{n}-u_{N_{x}+1,j}^{n+1}-u_{N_{x}+1,j}^{n}\right)\nonumber,
\end{eqnarray}
and analogously for the condition (\ref{BCtransmission_1Dhyperbolic}). 

\begin{remark}\label{remark_kk}
Note that, in order to ensure the positivity of the quantities in the above formulas deriving from the KK conditions, i.e. (\ref{BCtransmission2Dparabolicsourcetermg}) for the 2D domain and (\ref{BCparabolic1DKKsourceg}) or (\ref{impl_transm_hyp}) for the 1D domain, we also need to take care of the ratio between the KK coefficient $K$ and the space discretization steps. In particular, for (\ref{BCtransmission2Dparabolicsourcetermg}) and (\ref{impl_transm_hyp}) one needs to ensure that $K\frac{\triangle t}{\triangle x}$ and, respectively, $K\frac{\triangle t}{\triangle x}\triangle y$ is not too big in order to damp possible high obscillations produced by the term in parenthesis. Similarly, in (\ref{BCparabolic1DKKsourceg}) we need to check that $K\frac{\triangle t}{\triangle x} \sigma$ is small in order to prevent the growing of negative term.\\
Moreover, as previously discussed, we need to check that the numerical monotonicity conditions are satisfied:
\begin{align}
\frac{k_{1}}{\left(k_{2}+\varphi^n_{i,j}\right)^{\gamma}} \vert \partial^n_{x,i,j}\varphi^n_{i,j} \vert\leq\sqrt{D_{M}}\\
 \frac{k_{\omega} \omega}{1+\omega^n_{i,j}} T^n_{i,j} \leq 1 .
\end{align}
in the computational domain in order to ensure non-negative solutions.
\end{remark}

For the sake of completeness, we underline that at each time step a non-linear equation system must be solved, for which Newton-Krylov-subspace methods \cite{Knoll} can be used which take advantage of the mostly sparse structure of the jacobian matrix.\\

\section{Numerical tests and results}\label{sec:tests}

This section is devoted to the presentation of the numerical tests and the parameters of the problem are reported in Table \ref{table:param1}. Our aim is to show the ability of the simulation algorithm based on the model (\ref{eqsystem})-(\ref{GA1D}) to reproduce the qualitative behavior of the two population sharing the same habitat as observed in the videos of laboratory experiments.\\ 
We remark that we decided to do numerical simulations of the chip geometry assigning the 1D-hyperbolic-parabolic  model on channels since it seems more realistic. However, a numerical test on the behavior of the model (\ref{eqsystem})-(\ref{eqsystem1D}), with the doubly-parabolic model on channels, is provided in the last Example \ref{test4}.
 
\begin{table*}[tbp]\label{table}
\centering
{
\begin{tabular}{|p{2cm}|c|c|c|c|} \hline
 Parameter  &  Description& Units & Value& Ref.\\\hline\hline
 $D_{M}$  & Diffusivity of cells &  $\mu m^2/s$ & $9 \times 10^2$& \cite{murray} \\
 $D_T$& Diffusivity of cells &  $\mu m^2/s$ & $5.6\times 10^{1}$&\cite{murray} \\
 $D_{\varphi}, D_{\omega}$  & Diffusivity of chemoattractants & $\mu m^2/s$ &  $2 \times 10^2$& \cite{murray} \\ 
 $\alpha_T$ & decay rate of drug release & $s^{-1}$ &0 &- \\
$\alpha_M$ & decay rate of drug release & $s^{-1}$ & 0& -\\
$K_T$ & decay rate of $T$ caused by drug& $s^{-1}$ &0&- \\
$K_M$ & decay rate of $M$ caused by drug& $s^{-1}$  & 0& -\\
$\alpha_\varphi$ & growth rate of $\varphi$ & $ s^{-1}/cell$ &$10^{-1}$ & \cite{curk} \\
$\beta_\varphi$ &consumption rate of $\varphi$ & $ s^{-1}$ &$10^{-4}$& \cite{curk} \\
 $\alpha_\omega$ & growth rate of $\omega$& $ s^{-1}/cell$ &$10^{-1}$& \cite{curk} \\
$\beta_\omega$ & consumption rate of $\omega$& $ s^{-1}$ & $10^{-4}$&\cite{curk} \\
$k_1$ & cellular drift velocity & $ M cm^2 s^{-1}$& $3.9 \cdot 10^{-9}$& \cite{murray}\\
$k_2$ & receptor dissociation constant & $ M$ &$5 \cdot 10^{-6} $& \cite{murray}\\
$k_\omega$ & killing efficiency of immune cells& $\mu m/s$ per cell & 1& -\\
$\gamma$ & exponent in chemotactic response $\chi$ &  & 2&  \cite{murray}\\
$L$ & length of the corridor & $\mu m$  & 500& \\
$L_x$ & horizontal size of the box & $\mu m$  & 100 & \\
$L_y$ & vertical size of the box & $\mu m$  & 1000& \\
\hline\hline
\end{tabular} }
\vspace{0.2 in}    

\caption{Parameters of the problem.}\label{table:param1}
\end{table*}

\begin{example}\label{ex:accuracy}
Before we numerically simulate the laboratory experiment with the algorithm, we conduct a simple numerical test in order to prove its accuracy. We assumed the following setting: a left squared chamber $\Omega_{l}$ with one corridor positioned in the middle and only one cell family with initial distribution $u(x,y,0)=5e^{-\frac{1}{2}\left((x-0.5)^{2}+(y-0.5)^{2}\right)}$. Since we do not have any analytical solution for this problem, we choose $dt$ and $dx=dy$ small enough to obtain reasonable error estimations. In this case we use $dt=10^{-4}$ and $dx=dy=5 \times 10^{-4}$ for the approximation $u_{e}$ at time $t=100$ and calculate the error as the quantity $\| u_{e}-u_{\text{approx}}\|$ in $L^1$-norm.\\
In order to confirm the order of our scheme, we use a log-log-plot with constant and small enough $dt$ (resp. $dx$), and decreasing $dx$ (resp. $dt$).
As shown in Figures \ref{fig:errorplotdt} and \ref{fig:errorplotdx} the time order and space order equals to line with slope 2 in the log-log plot which corresponds to our scheme of order 2 in space and time.\\
\begin{figure}[h!]
\centering
\includegraphics[scale=0.3]{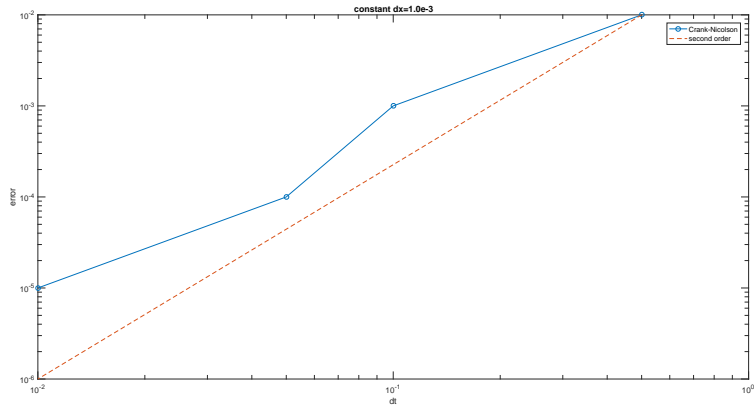}
\caption{Log-log plot of the error, namely the quantity $\| u_{e}-u_{\text{approx}}\|$ in $L^1$-norm as a
function of the space step, with fixed $dt=10^{-3}$ and decreasing $dx=0.5, 0.1, 0.05, 0.001$ at time $t=100$. We depict in blue the obtained error and in red a line with slope 2 for comparison.}
\label{fig:errorplotdt}
\end{figure}
\begin{figure}[h!]
\centering
\includegraphics[scale=0.3]{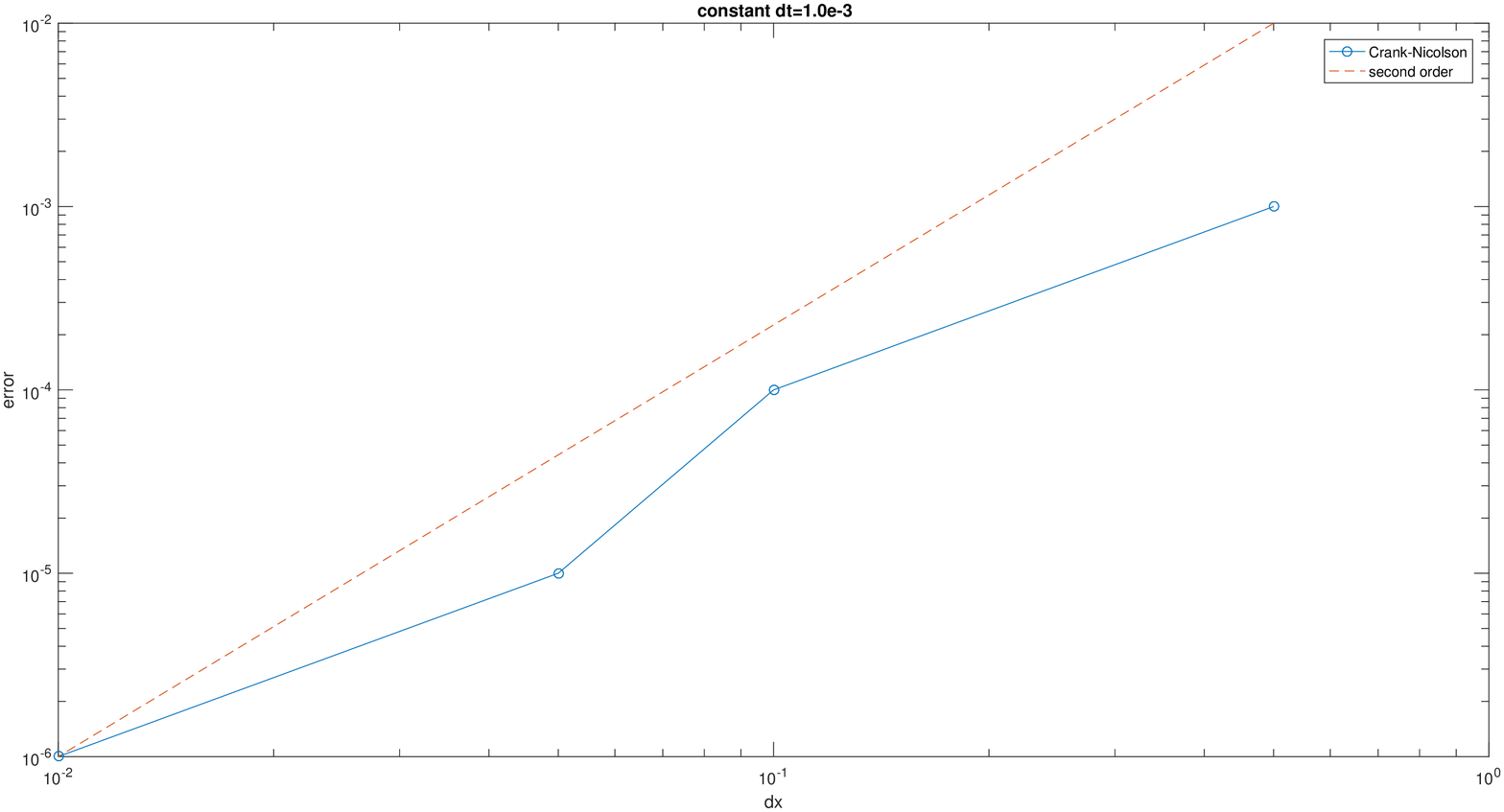}
\caption{Log-log plot of the error, namely the quantity $\| u_{e}-u_{\text{approx}}\|$ in $L^1$-norm as a
function of the time step, with fixed $dx=10^{-3}$ and decreasing $dt=0.5, 0.1, 0.05, 0.001$ at time  $t=100$. We depict in blue the obtained error and in red a line with slope 2 for comparison.}
\label{fig:errorplotdx}
\end{figure}
\end{example}
Now we describe the simulation of the chip environment. All the simulations were performed in MATLAB\textcircled{c}. The computational time for a simulation on the complete geometry until time $t=100$, takes about 40 seconds on an Intel(R) Core(TM) i7-3630 QM CPU 2.4 GHz.

\begin{example}\label{test1}
For the following numerical simulation we replicate the laboratory experiment by having the two domains (representing the two chambers) $\Omega_{l}=[0,L_{x}]\times[0,L_{y}]$ and $\Omega_{r}:=[L_{x}+L,2L_{x}+L]\times[0,L_{y}]$ and 5 corridors $I_{m}:=[0,L], \ m=1,\ldots,5$  having the same width $\sigma$ and equispaced from each other.\\
We choose $L_{x}=100\mu m$, $L_{y}=1000\mu m$, $L=500\mu m$ and $\sigma=12\mu m$ to accurately simulate the experiment shown on the video footage.\\
The initial condition (time $t=0$) for the tumor cells distribution on the chip for $(x,y) \in \Omega_l$ is chosen as:
\begin{equation}
T(x,y,0)=5e^{-\frac{1}{2}\left( x^{2}+y^{2}\right)}+5e^{-\frac{1}{2}\left( x^{2}+(y-5)^{2}\right)}+5e^{-\frac{1}{2}\left( x^{2}+(y-10)^{2}\right)},
\end{equation}
whereas in the corridors and the right chamber no tumor cells are present.\\
For the immune cells distribution on the chip for $(x,y) \in \Omega_r$ we assign:
\begin{equation}
M(x,y,0)=5e^{-\frac{1}{2}\left((x-1-L_{x}-L)^{2}+(y-5)^{2}\right)},
\end{equation}
 whereas no immune cells are present in the left chamber nor in the corridors.\\
 For the chemoattractants we set a constant initial density: $\omega(x,y,0)=0$ (all domains) and $\varphi(x,y,0)=2$ for $x,y \in \Omega_{l}\setminus \{L_x \times [a_2,b_2]\}$, $\varphi(x,y,0)=0$ for $x,y \in \Omega_{r}-[a_2,b_2]$ and a linear decreasing in space initial value $\varphi(x,y,0)=-0.01 x+5$ in correspondence of nodes $2L-2R$, namely at $(L_x,y)$ and $(L_x+L,y)$ for $y \in [a_2,b_2]$.\\
For this simulation test we choose the parameters for each domain as given in Table \ref{table:param1} for both the chambers and for all the corridors we used the same parameters. \\
The numerical method implemented is listed in paragraph \ref{sec:algo}; for the 1D channels the AHO-Scheme (\ref{AHO2_scheme_cells_impl}) is implemented since we are considering the hyperbolic-parabolic model. The discretization grid has time step size $\triangle t=10^{-3}$ and space size $\triangle x=\triangle y=0.25$.\\
\begin{figure}[H]
\centering
\includegraphics[scale=0.75]{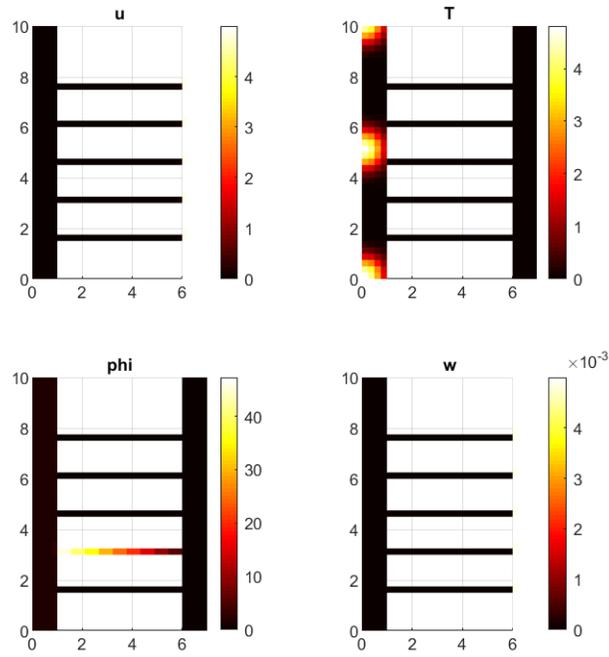}
\caption{Initial distribution of tumor cell densities $T$ and of immune cells $u$ at time $t=0$.}
\label{fig:plott=0}
\end{figure}
\begin{figure}[H]
\centering
\includegraphics[scale=0.75]{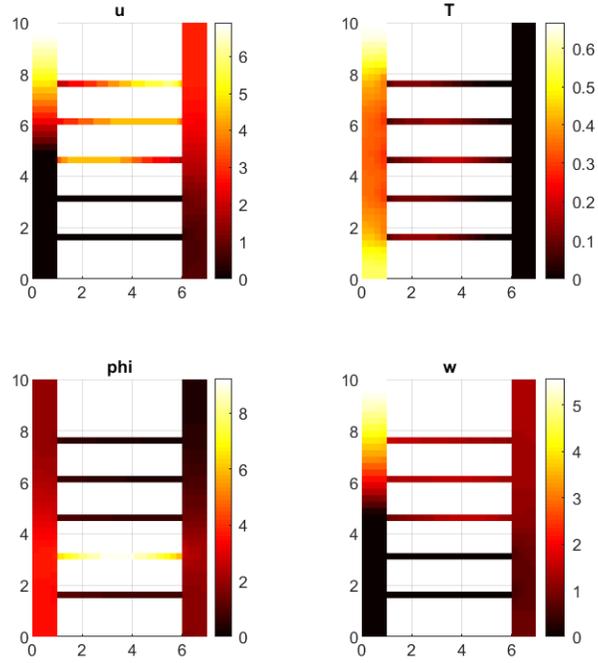}
\caption{Simulation of model (\ref{eqsystem})-(\ref{GA1D}). On the left: tumor cell densities $T$. On the right: immune cells $u$ diffusing around right chambers and entering corridors in higher quantities than $T$ at time $t=5$.}
\label{fig:plott=5}
\end{figure}
\begin{figure}[H]
\centering
\includegraphics[scale=0.75]{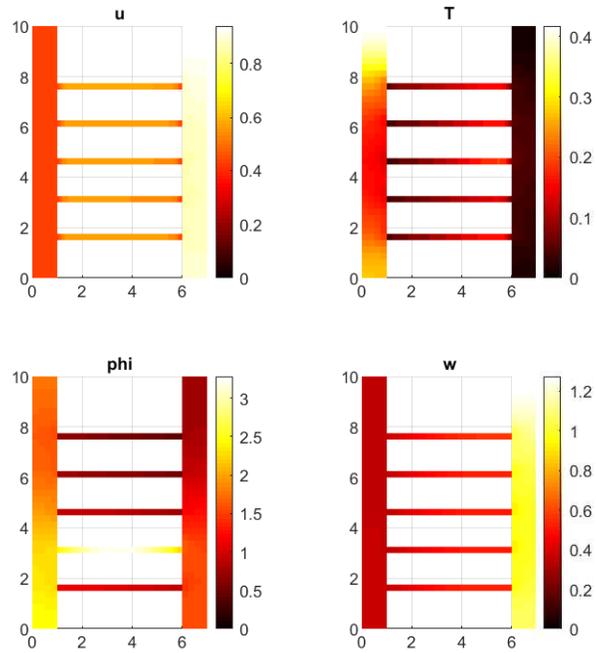}
\caption{Simulation of model (\ref{eqsystem})-(\ref{GA1D}). On the left: tumor cells $T$ get pushed towards the top of the left chamber. On the right: immune cells accumulating in the top of right chamber at time $t=10$.}
\label{fig:plott=10}
\end{figure}
\begin{figure}[H]
\centering
\includegraphics[scale=0.75]{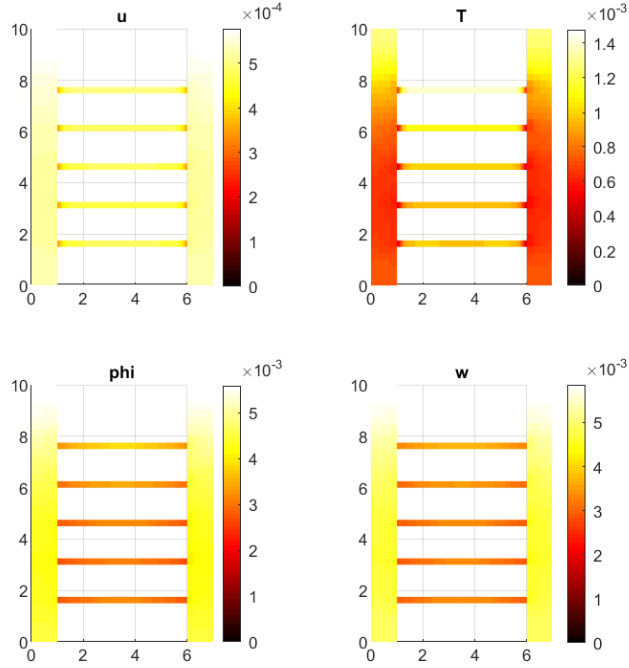}
\caption{Simulation of model (\ref{eqsystem})-(\ref{GA1D}). On the left: tumor cells $T$ reach densities values of order $10^{-5}$ but accumulate mostly around top of both chambers. On the right: distribution of immune cells $u$ at time $t=100$. }
\label{fig:plott=100}
\end{figure}

In Figures \ref{fig:plott=0}, \ref{fig:plott=5}, \ref{fig:plott=10} and \ref{fig:plott=100} we can see the density of the tumor cells $T$ and immune cells $u$ for different times $t=0$, $t=5$, $t=10$ and $t=100$ accordingly. Note that at time $t=0$ tumor cells are present in the left chamber only and immune cells are present in the right chamber only.\\
Since no chemoattractant is present at the initial time $t=0$, both cells diffuse around their chamber and slowly entering the corridors while creating chemoattractant $\varphi$ and $\omega$.\\
But already at time $t=5$ we notice that in the middle of the left chamber the tumor cells are getting pushed towards the buttom and top of their chambers. Indeed, tumor cell densities $T$ slowly diffuse around the left chamber but accumulate around the bottom and top and partly enter the corridors.
 This is due to the fact that the chemoattractant $\varphi$ produced by cancer cells (the annexin) induce a migration of the immune cells $M$ towards the tumor cells $T$ causing a higher migration towards the center of the left chamber where the initial distribution of tumor cells was closest to the chambers.\\
For $t=10$ in Fig. \ref{fig:plott=10} and $t=100$ in Fig. \ref{fig:plott=100} we see that, since most tumor cells are accumulating on the top, they manage to diffuse through the nearest corridor on the top into the right chamber; on the other hand, the immune cells continue to migrate towards the highest concentration of chemoattractant $\varphi$, which is where the tumor cells $T$ concentrate.\\
Especially at time $t=100$, we can see that the quantity of tumor cells has dramatically decreased compared to the immune cells caused by the action of chemokine $\omega$.
\end{example}
\begin{example}\label{test2}
In this numerical test we used the same settings of Example \ref{test1}, where the only difference consists in a much stronger chemotaxis, i.e. $k_1$  in $\chi(M,\varphi)$ is 50 times larger than in the previous Example \ref{test1}.
The results are depicted in the following Figures \ref{fig:plott1=0}-\ref{fig:plott1=100}.
\begin{figure}[H]
\centering
\includegraphics[scale=0.75]{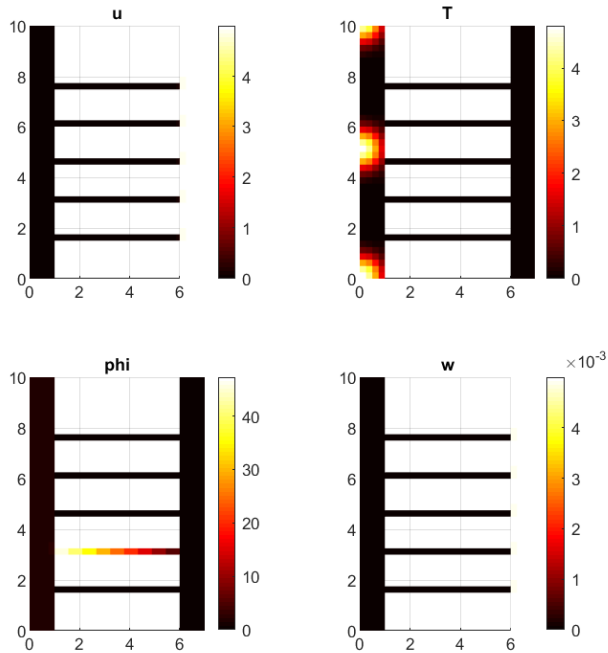}
\caption{Initial distribution of tumor cell densities $T$ and of immune cells $u$ at time $t=0$ for the model with stronger chemotaxis.}
\label{fig:plott1=0}
\end{figure}
\begin{figure}[H]
\centering
\includegraphics[scale=0.75]{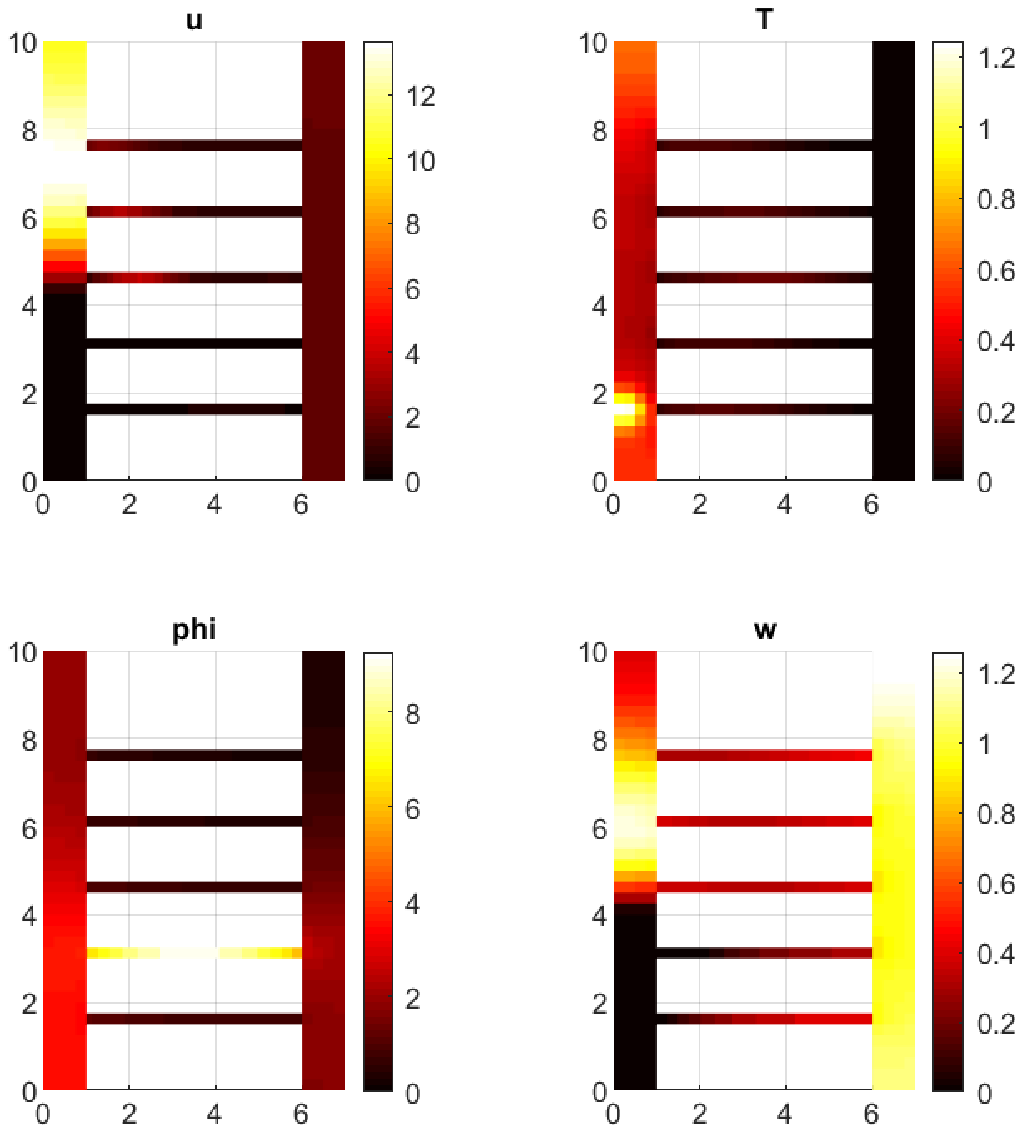}
\caption{Densities time $t=5$ for model (\ref{eqsystem})-(\ref{GA1D}) with stronger chemotaxis.}
\label{fig:plott1=5}
\end{figure}
\begin{figure}[H]
\centering
\includegraphics[scale=0.75]{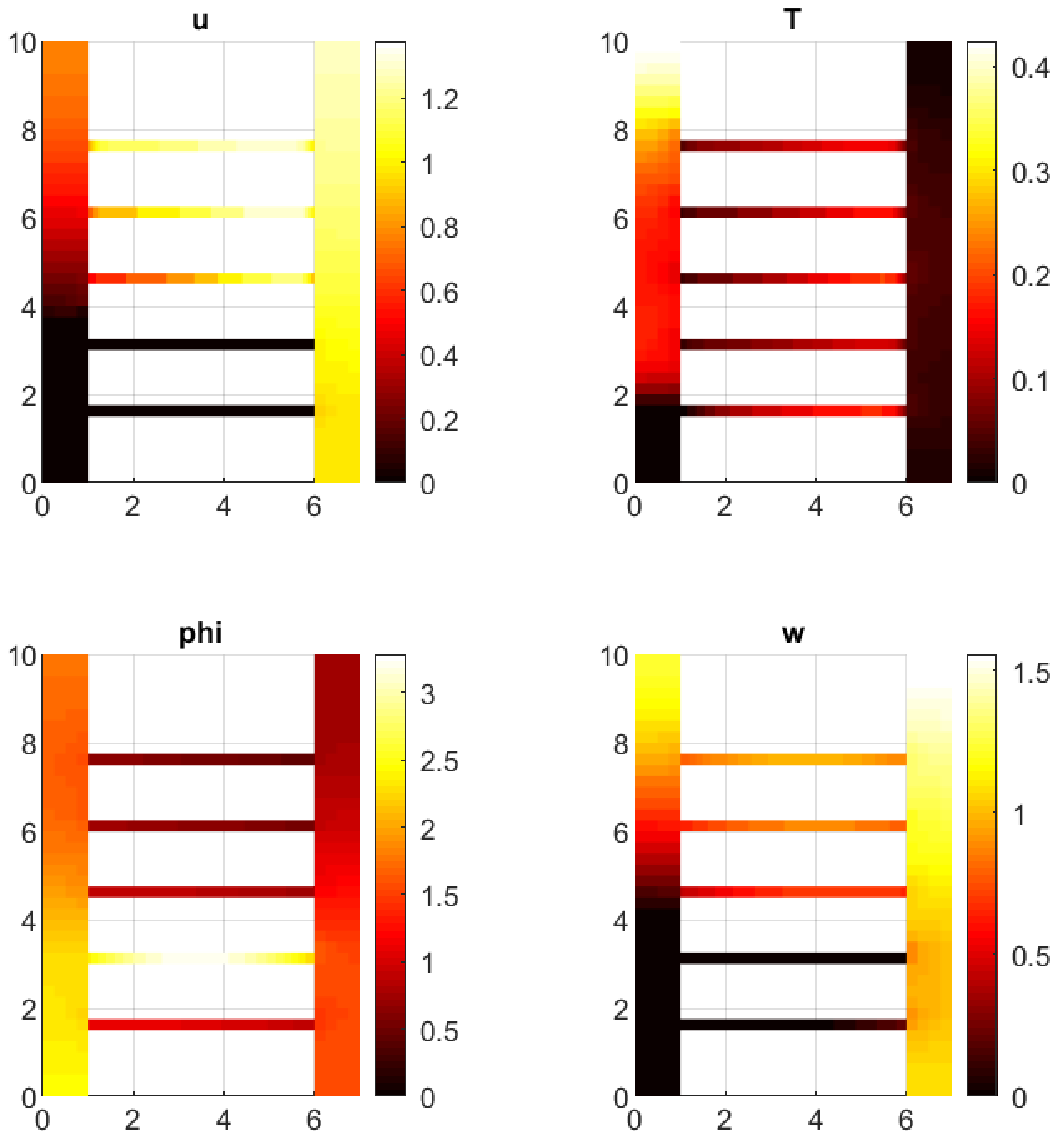}
\caption{Densities at time $t=10$ for model (\ref{eqsystem})-(\ref{GA1D}) with stronger chemotaxis.}
\label{fig:plott1=10}
\end{figure}
\begin{figure}[H]
\centering
\includegraphics[scale=0.75]{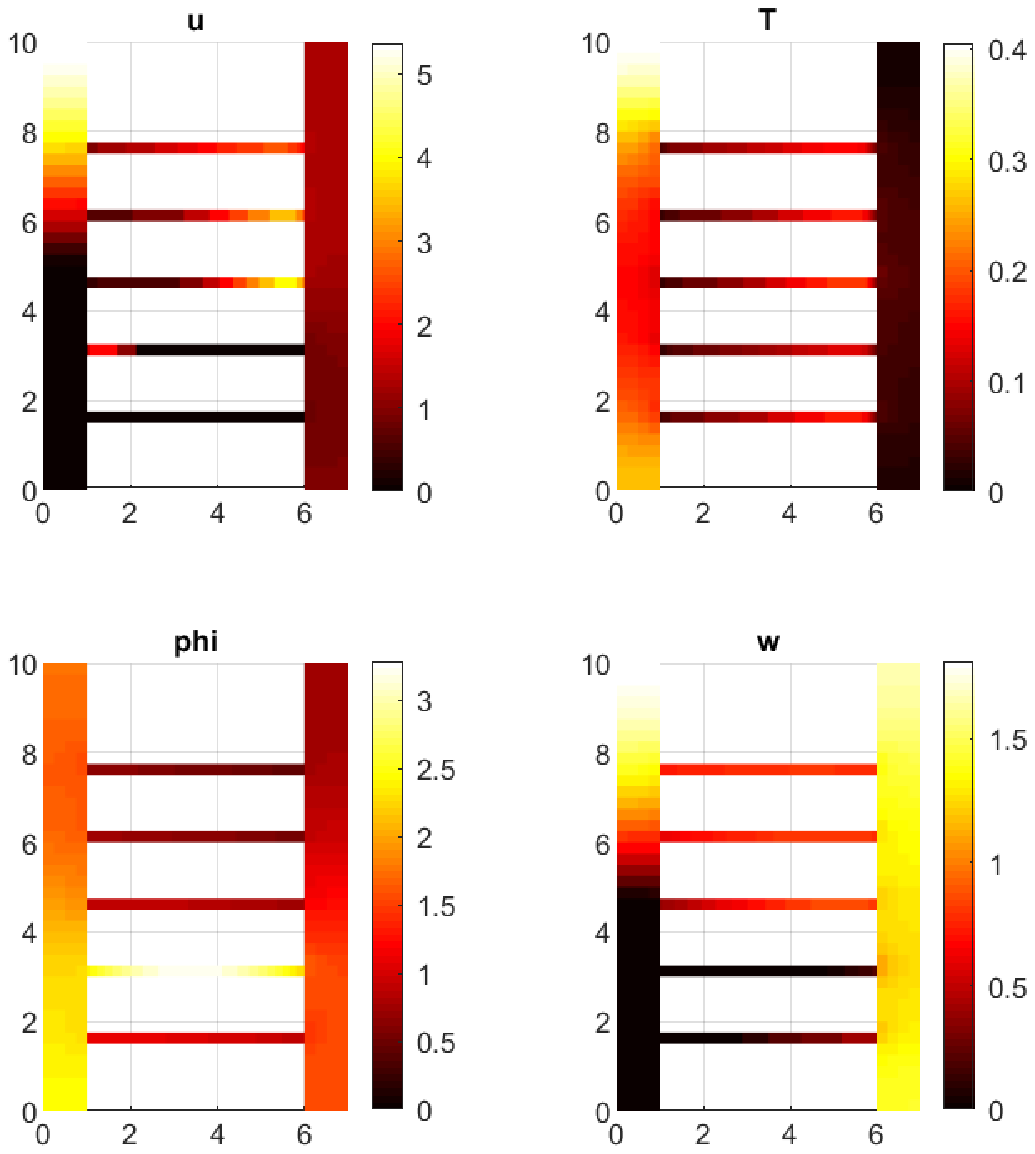}
\caption{Densities at time $t=100$ for model (\ref{eqsystem})-(\ref{GA1D}) with stronger chemotaxis.}
\label{fig:plott1=100}
\end{figure}
\begin{figure}[H]
\centering
\includegraphics[scale=0.75]{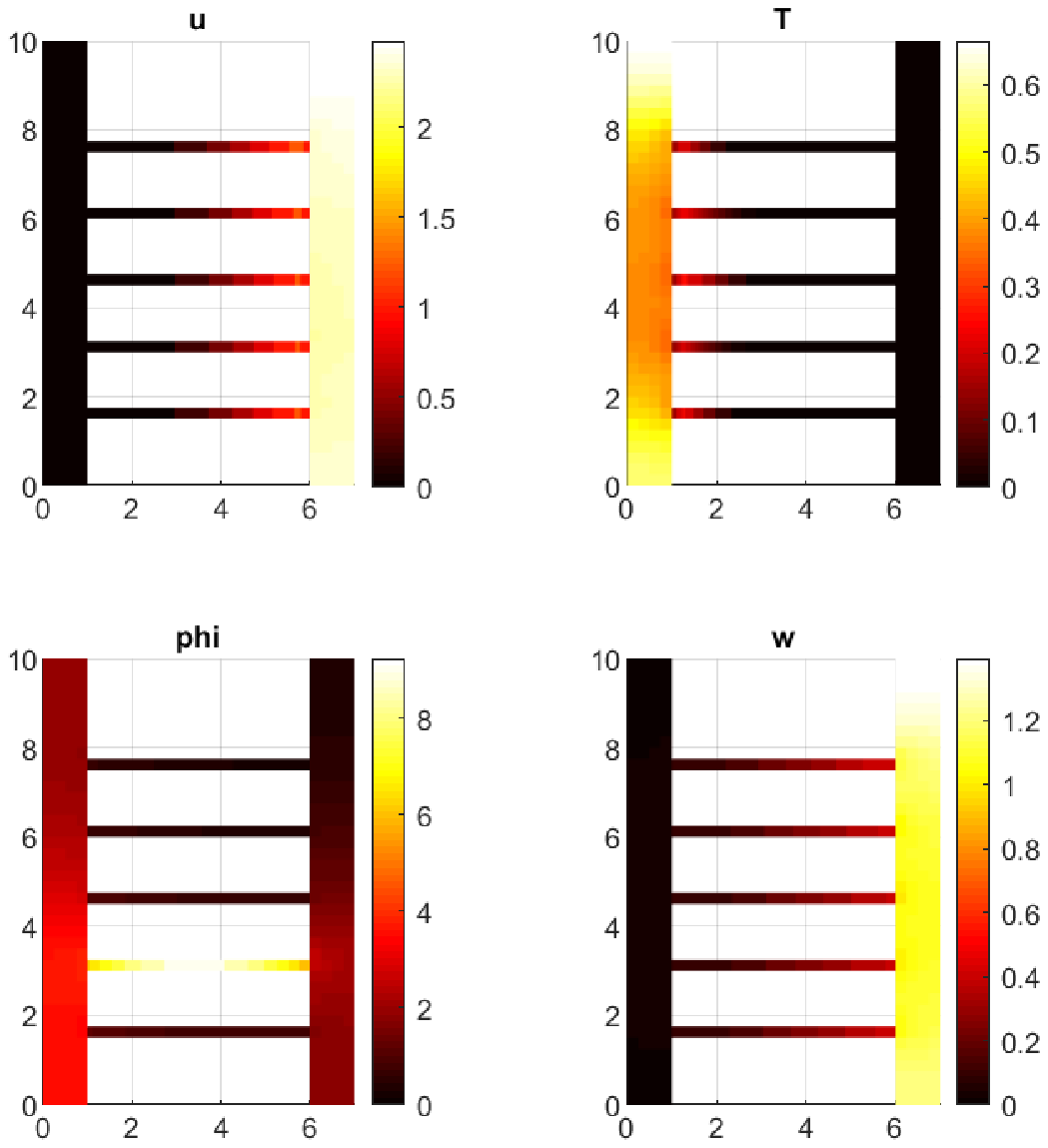}
\caption{Densities time $t=5$ for the 2D-1D-doubly-parabolic model (\ref{eqsystem})-(\ref{eqsystem1D}).}
\label{fig:plott2=5}
\end{figure}
\begin{figure}[H]
\centering
\includegraphics[scale=0.75]{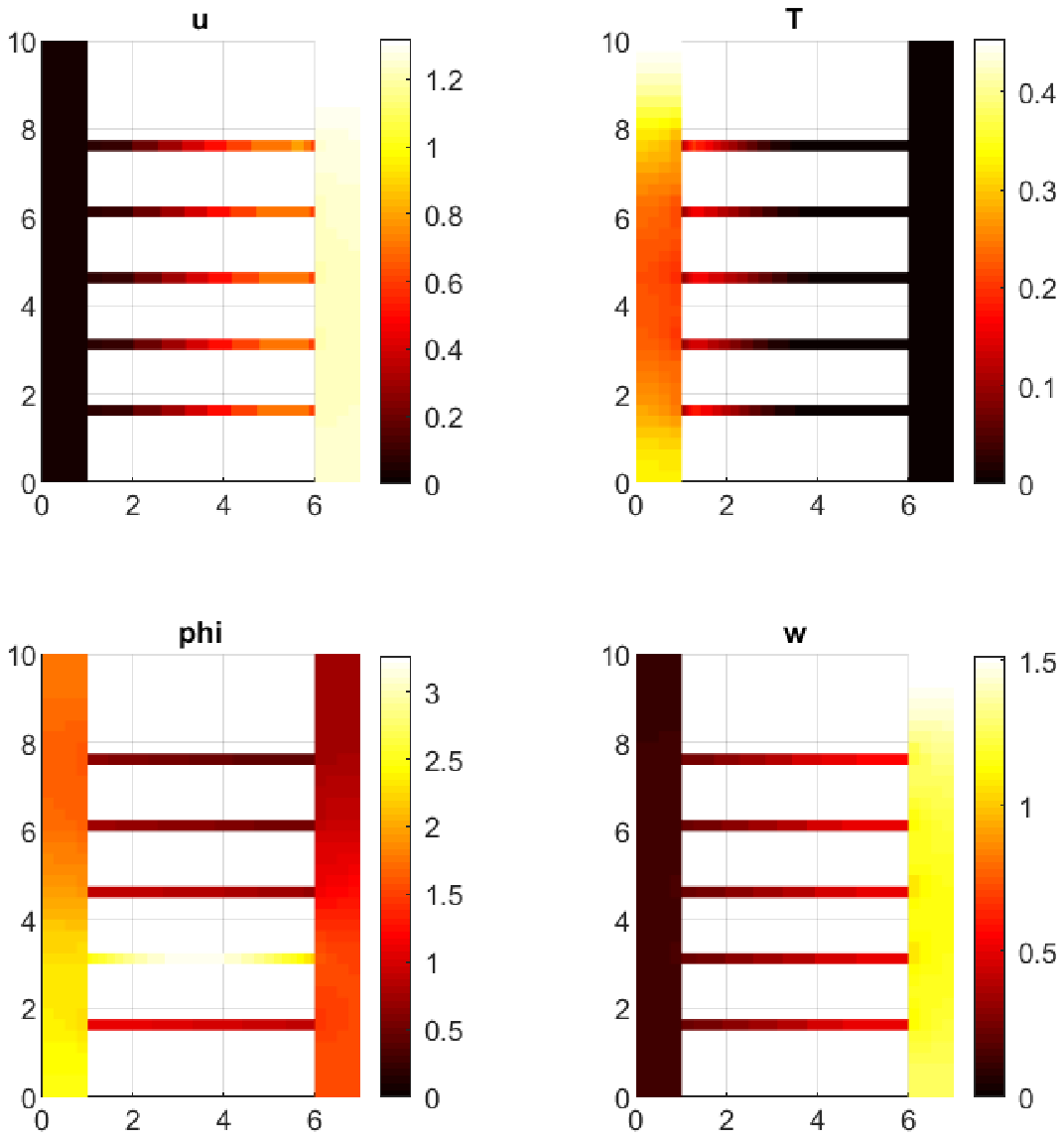}
\caption{Densities at time $t=10$ for the 2D-1D-doubly-parabolic model (\ref{eqsystem})-(\ref{eqsystem1D}).}
\label{fig:plott2=10}
\end{figure}
\begin{figure}[H]
\centering
\includegraphics[scale=0.75]{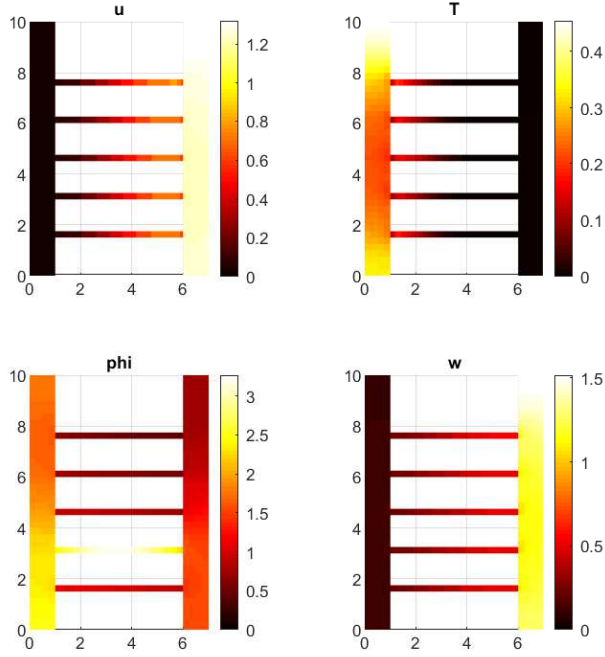}
\caption{Densities at time $t=100$ for the 2D-1D-doubly-parabolic model (\ref{eqsystem})-(\ref{eqsystem1D}). }
\label{fig:plott2=100}
\end{figure}
We can clearly see in Fig. \ref{fig:plott1=5} compared to Fig. \ref{fig:plott=5}, that the immune cells $u$ are much more massively moving towards the left chamber due to the chemoattractant $\varphi$, which causes a slightly higher concentration in the left chamber. But due to the diffusion of the chemoattractant $\varphi$ and the creation of more chemoattractant from the tumor cells, the graphs of both Examples \ref{test1} and \ref{test2} are getting similar during the time evolution until a difference is no more noticeable, around time $t=100$.
\end{example}

\begin{example}\label{test4}
In this last Example, we tested the 1D-doubly-parabolic model on channels and compared it with the hyperbolic-parabolic model used in the previous Examples.
In Figures \ref{fig:plott2=5}-\ref{fig:plott2=100} we assigned the 2D-doubly-parabolic model which uses the parabolic partial differential equation to describe the movement in the 1D channels.\\
By using the same initial data as for the other Examples, we notice that for time $t=100$, the doubly-parabolic model in Fig. \ref{fig:plott2=100} seems to have a similar pattern as for the hyperbolic-parabolic model depicted in Fig. \ref{fig:plott=100} and the hyperbolic-parabolic model with stronger chemotaxis in Fig. \ref{fig:plott1=100}, but the scale differs a lot between these models.\\
Whereas we have for the tumor cells $T$ a maximum concentration of $10^{-5}$ for the hyperbolic-parabolic model, and $10^{-7}$ for the hyperbolic-parabolic mode with stronger chemotaxis, we see clearly that for the doubly-parabolic model, the concentration of the tumor cells $T$ is of the order of $10^{-2}$. This is due to the much slower movement of the immune cells through the corridors. This also explains the much higher concentration of the chemoattractant $\phi$ because of the much higher concentration of $T$ compared to the other models.
\end{example}

In the following Figure \ref{fig:cellvisualization}, we represent the density of tumor cells and immune cells depicted in Figures \ref{fig:plott=0}-\ref{fig:plott=100} as individuals, by randomly placing them according to their density. The higher the density at a given point, the more cells will be distributed randomly around that area. If the density is lower than a chosen threshold in a certain point, no cells will be represented around it. 

\begin{figure}[h!] 
    \centering
    \begin{subfigure}[t]{0.5\textwidth}
        \centering
        \includegraphics[height=1.2in]{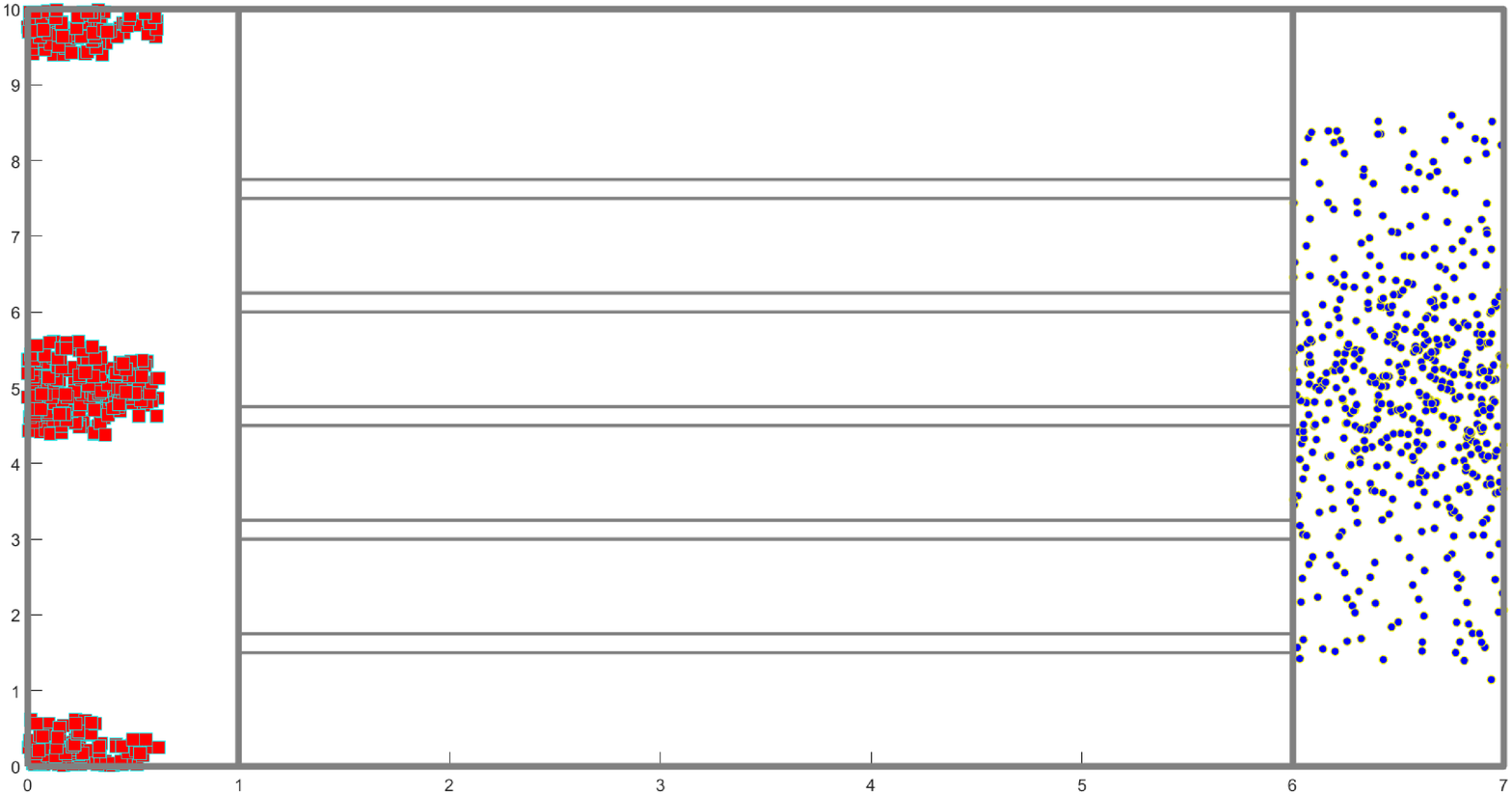}
        \caption{Visualization for time t=0.}
    \end{subfigure}
    \begin{subfigure}[t]{0.5\textwidth}
        \centering
        \includegraphics[height=1.2in]{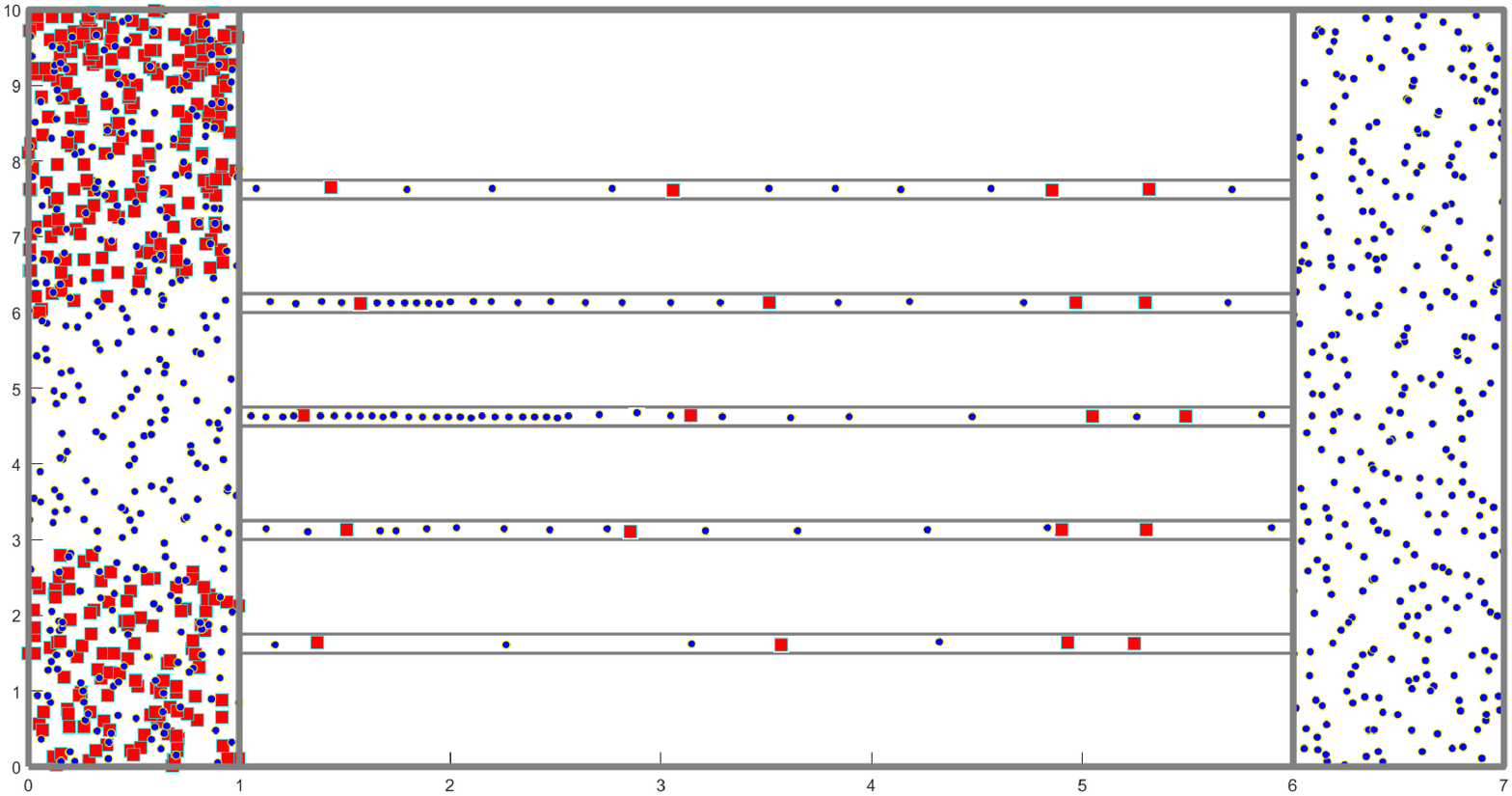}
        \caption{Visualization for time t=5.}
    \end{subfigure}
    ~
     \begin{subfigure}[t]{0.5\textwidth}
        \centering
        \includegraphics[height=1.2in]{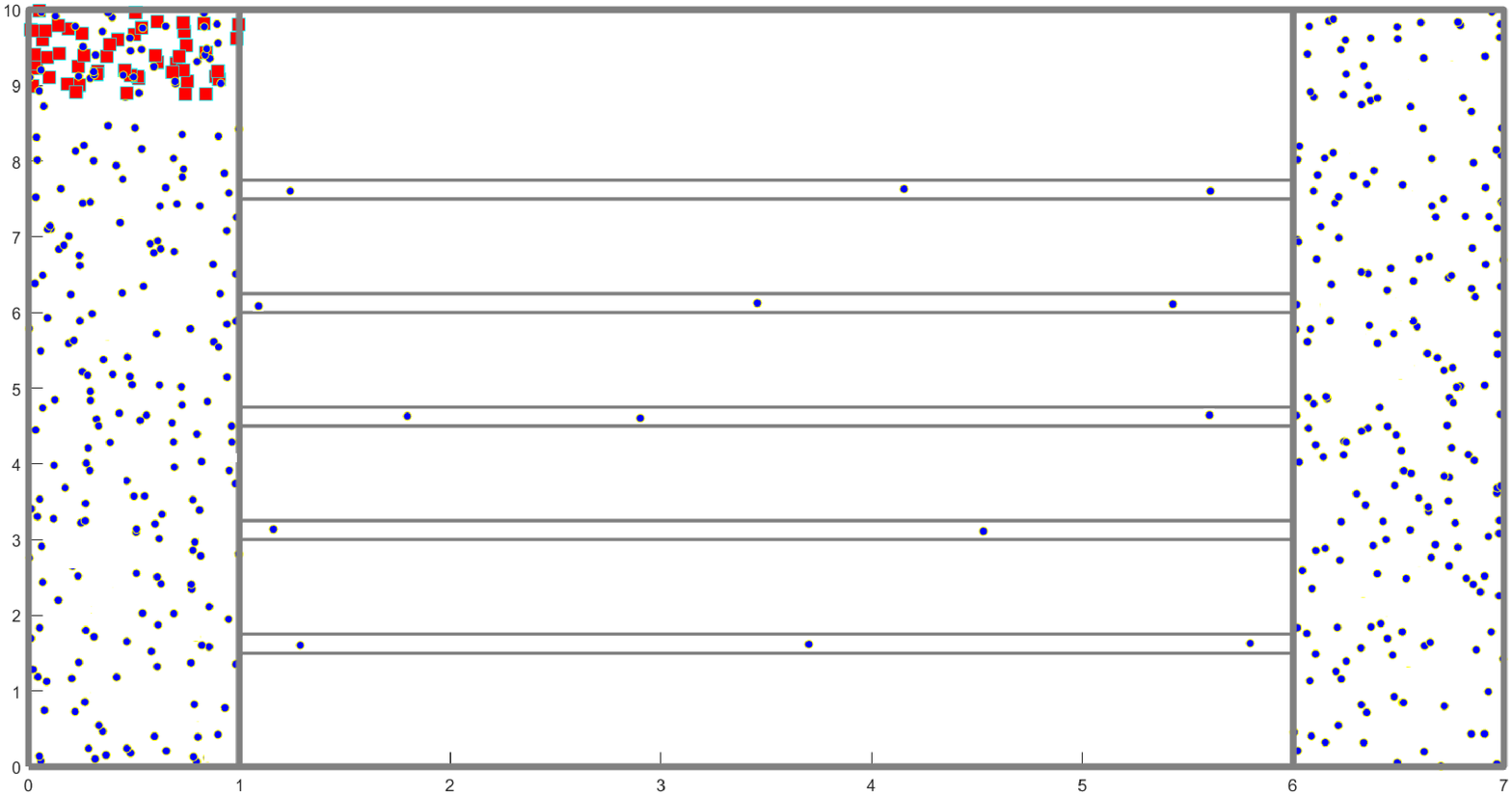}
        \caption{Visualization for time t=50.}
    \end{subfigure}
    \caption{Visualization of immune cells (blue dots) and tumor cells(red squares) for time t=0, t=5 and t=50 by using the density of each quantity and representing them as cells.}\label{fig:cellvisualization}
\end{figure}

\section{Conclusion and future perspectives}\label{sec:concl}
The principal feature of the present work has been the development of a simulation tool to describe cell movements and interactions inside microfluidic chip environment. Our study focused on both the modelling and the numerical point of view. Indeed, schematizing the chip geometry as two 2D-boxes connected by a network of 1D-channels, the main issues were:
\begin{itemize}
\item the introduction of mass-preserving conditions involving the balancing of incoming and outgoing fluxes passing through interfaces between 2D and 1D domains;
\item the development of mass-preserving numerical schemes at the boundaries of 2D domain and mass-preserving transmission conditions at the 2D-1D interfaces.  
\end{itemize}
Furthermore, from the modelling point of view, we studied the dynamics in the channels in case of doubly-parabolic model and hyperbolic-parabolic model. Since we obtained comparable asymptotic states, we decided to apply the hyperbolic-parabolic model in order to have finite speed of propagation in the channels which seems to be more realistic.
In this framework, having in mind the laboratory experiments on chip described in section \ref{sec:bio}, it was possible to simulate the chip environment with two species of living cell moving in it. Moreover, we remark that we can simulate more complicated situations where more than two cell species are present.\\
As a further development of the present study, we will work on the calibration of the model against experimental data.


\end{document}